\documentclass[11pt]{article} \usepackage{amsmath}
\usepackage{amssymb} \usepackage{amscd} \usepackage{amsthm}

\def\endproof{\relax\ifmmode\expandafter\endproofmath\else
  \unskip\nobreak\hfil\penalty50\hskip.75em\hbox{}\nobreak\hfil\bull
  {\parfillskip=0pt \finalhyphendemerits=0 \bigbreak}\fi} 
\def\endproofmath$${\eqno\bull$$\bigbreak}
\def\bull{\vbox{\hrule\hbox{\vrule\kern3pt\vbox{\kern6pt}\kern3pt\vrule}\hrule}}
\newtheorem{theorem}{Theorem}[subsection]
\newtheorem*{main}{Theorem} 
\newtheorem{proposition}[theorem]{Proposition}

\newtheorem{lemma}[theorem]{Lemma}
\newtheorem{claim}[theorem]{Claim}

\newtheorem{corollary}[theorem]{Corollary}

\newtheorem{D}[theorem]{Definition}
\newenvironment{defn}{\begin{D} \rm }{\end{D}}

\newtheorem{R}[theorem]{Remark}
\newenvironment{remark}{\begin{R}\rm }{\end{R}}

\def\Zee{\mathbb{Z}}

\def\Cee{\mathbb{C}} 
\def\Aff{\mathbb{A}}
\def\Pee{\mathbb{P}} 
 
\def\cS{\mathcal{S}}
\def\Gm{\mathbb{G}_m}
\def\scrO{\mathcal{O}} 
 
\def\spcheck{^{\vee}}
\def\frak{\mathfrak} 
\def\Pic{\operatorname{Pic}}
\def\Ker{\operatorname{Ker}}

\def\Hom{\operatorname{Hom}}
 
\def\Aut{\operatorname{Aut}}
 
\def\Id{\operatorname{Id}}
 
\def\Spec{\operatorname{Spec}}
\def\ad{\operatorname{ad}} 
\def\Ad{\operatorname{Ad}}
\def\Lie{\operatorname{Lie}}

\title{ Automorphism Sheaves, Spectral Covers, and the Kostant and Steinberg Sections}
\author{Robert Friedman\thanks{The first author was partially
    supported by NSF grant DMS-02-00810.}
\  and John W. Morgan\thanks{The second author was partially supported
   by NSF grant DMS-01-03877.}}

\begin{document}
\maketitle

\section*{Introduction}

Throughout this paper, $G$ denotes a simple and simply connected algebraic
group over $\Cee$ of rank $r$ and $H$ is a Cartan subgroup, with Lie algebras $\frak
g=\Lie G$ and
$\frak h= \Lie H$. Let $R$ be the root system of the pair $(G,H)$, $W$ the Weyl group,
and $\Lambda\subseteq \frak h$ the coroot lattice. Fix once and for all a positive
Weyl chamber, i.e.\ a set of simple roots $\Delta$. The geometric invariant theory
quotient of
$\frak g$ by the adjoint action of
$G$ is identified with
$\frak h/W$, which by a theorem of Chevalley is isomorphic to $\Aff^r$. Let $\frak
g_{\rm reg}\subseteq \frak g$ be the Zariski open and dense subset of $\frak g$
consisting of regular elements, i.e.\ of elements $X$ whose centralizer (in $\frak g$
or $G$) has the minimal dimension $r$. If $X$ is regular, then the centralizers of $X$
in $\frak g$ and in $G$, i.e.\ the stabilizer of $X$ in $G$ under the adjoint
representation, are abelian. One can show that each fiber of the adjoint quotient
morphism
$\frak g\to \frak h/W$ contains a unique orbit of regular elements, so that the
restriction of the adjoint quotient morphism to $\frak g_{\rm reg}$ induces an
isomorphism $\frak g_{\rm reg}/G
\to \frak h/W$. In 1963, Kostant \cite{Kostant} proved that the morphism
$\frak g_{\rm reg} \to \frak h/W$ is smooth and constructed a section
$\sigma \colon \frak h/W \to
\frak g_{\rm reg}$ of the adjoint quotient morphism, generalizing the existence of
rational canonical forms of a matrix. Actually, the construction yields a family of
such sections, defined as follows. Let $X \in\frak g_{\rm reg}$ be a principal
nilpotent matrix and let $L\subseteq \frak g$ be a linear complement to
$\operatorname{Im}\ad X$. Then the affine subvariety $L+X$ is contained in $\frak
g_{\rm reg}$ and maps isomorphically to $\frak h/W$, and so defines a section of the
adjoint quotient morphism.

One can make an analogous construction for the adjoint quotient of $G$, i.e.\ the
geometric invariant theory quotient of the action of $G$ on itself by conjugation. In
this case, the adjoint quotient is equal to $H/W$, which is again isomorphic to
$\Aff^r$ \cite{Bour}. Let $G_{\rm reg}$ be the Zariski open and dense subset of regular
elements, where by definition $g\in G$ is regular if and only if the dimension of its
centralizer is $r$, the minimal possible dimension. If $g$ is regular, then its
centralizer is abelian. Then again each fiber of the morphism $G\to H/W$ contains a
unique conjugacy class of regular elements, so that $G_{\rm reg}/G$ is isomorphic to
$H/W$. In 1965, using methods very different from those of Kostant, Steinberg
\cite{Steinberg} proved that the morphism
$G_{\rm reg} \to H/W$ is smooth and constructed a section
$\sigma \colon H/W \to G_{\rm reg}$ of the adjoint quotient morphism. Although the
constructions are quite different, they are related as follows. Let $\sigma \colon H/W
\to G_{\rm reg}$ be any section of the adjoint quotient morphism, and let $X\in \frak
g_{\rm reg}$ be the unique principal nilpotent matrix such that $\exp X = \sigma(e)$,
where we denote by $e$ the image of the identity in $H/W$. Then
$\operatorname{Im}d\sigma_e$ is a linear subspace $L$ of $\frak g$ which is a
complement to $\operatorname{Im}\ad X$, and so $L+X$ defines a Kostant section of the
adjoint quotient morphism for $\frak g$.

In \cite{FMIII}, the authors gave another method for constructing sections of the
adjoint quotient morphism, either for $\frak g$ or for $G$, as well as
generalizations to other contexts, which will be explained below. It is natural, then,
to ask if these sections are essentially the same as those defined by Kostant and
Steinberg, and more generally whether the sections of Kostant and Steinberg are in some
sense unique. It follows from the work of Kostant and Steinberg that any such section
must take values in the open set of regular elements. Clearly, if
$\sigma \colon \frak h/W \to \frak g_{\rm reg}$ is such a section and $f\colon \frak
h/W \to G$ is any morphism, then $\ad f(\sigma)\colon \frak h/W \to \frak g_{\rm reg}$
is again a section, and similarly for sections of $G\to H/W$. We shall refer to two
sections which are related in this way as \textsl{conjugate}. The main result of this
paper is then the following:

\begin{main} All sections of the adjoint quotient morphism $\frak g \to \frak
h/W$ are conjugate. Similarly, all sections of the adjoint quotient morphism $G \to
H/W$ are conjugate.
\end{main}

These problems can be translated into a more geometric context. Let $E$ be an
irreducible plane cubic curve with a cusp singularity and a marked point  $p_0 \in
E_{\rm reg}$. Then there is a natural group structure on $E_{\rm reg}$ with origin
$p_0$, and it is isomorphic to the additive group $\mathbb{G}_a$. Let $n\colon
\widetilde E\cong \Pee^1\to E$ be the normalization of $E$. Let $(\xi, \theta)$ be a
pair consisting of a principal $G$-bundle on $E$ and an isomorphism of principal
bundles $\theta\colon n^*\xi \to \widetilde E \times G$. By \cite[Corollary
3.1.2]{FMIII}, such pairs are classified by elements $X\in \frak g$, and hence
principal $G$-bundles $\xi$ such that $n^*\xi$ is trivial are classified by $G$-orbits
in $\frak g$. If the bundle $\xi$ corresponds to the $G$-orbit containing $X$, then the
group  $\Aut \xi$ of automorphisms of the principal bundle $\xi$ is isomorphic to the
stabilizer of $X$ in $G$. Thus, if $\xi$ is a $G$-bundle over $E$ such that
 $n^*\xi$ is trivial,  we define  $\xi$ to be
\textsl{regular} if  $\Aut \xi$ has dimension $r$. More
generally, let $B$ be a scheme over $\Cee$ and let $(\Xi, \Theta)$ consist of a
principal $G$-bundle over $B\times E$ and a bundle isomorphism $\Theta \colon (\Id\times
n)^*\Xi \cong B\times \widetilde{E}\times G$. Then such pairs are classified by
morphisms
$\varphi\colon B \to \frak g$. Given two such pairs $(\Xi, \Theta)$ and $(\Xi',
\Theta')$ corresponding to morphisms $\varphi$ and $\varphi'$ from $B$ to $\frak g$,
the bundles $\Xi$ and $\Xi'$ are isomorphic if and only if there exists a morphism
$f\colon B \to G$ such that $\ad f(\varphi) =\varphi'$. Applying this construction in
particular to a Kostant section $\sigma \colon \frak h/W \to \frak g_{\rm reg}$ defines
a $G$-bundle $\Xi \to (\frak h/W)\times E$ whose restriction to each slice
$\{x\}\times E$ is regular. Clearly, $\frak h/W$ is a coarse moduli space for the
set of bundles which become trivial on $\widetilde E$, $\Xi$ is a ``universal" bundle in
an appropriate sense, and our theorem becomes the statement that there is a unique such
bundle up to isomorphism.

There is an analogous construction in the case of the group $G$. In this case, let $E$
be an irreducible plane curve with a node and a marked point  $p_0 \in
E_{\rm reg}$. Then the  natural group structure on $E_{\rm reg}$ with origin
$p_0$  is isomorphic to   $\Gm$. Pairs $(\xi, \theta)$  consisting of a
principal $G$-bundle on $E$ and an isomorphism of principal bundles $\theta\colon
n^*\xi \to \widetilde E \times G$ are classified by an element $g\in G$, and more
generally pairs $(\Xi, \Theta)$ consisting of a
principal $G$-bundle over $B\times E$ and a bundle isomorphism $\Theta \colon (\Id\times
n)^*\Xi \cong B\times \widetilde E\times G$  are classified by morphisms
$\varphi\colon B \to G$. Regularity is defined as before. Using the Steinberg section
now to define a $G$-bundle over $(H/W)\times E$, 
$H/W$ is a coarse moduli space for the set of bundles which become trivial on
$\widetilde E$,
$\Xi$ is a ``universal" bundle in an appropriate sense, and our theorem is again the
statement that there is a unique such bundle up to isomorphism.

These constructions can be further generalized. Aside from nodal and cuspidal cubics,
one can also consider smooth cubic curves $E$ together with an origin $p_0$, i.e.\
elliptic curves. The group law on  $E$ gives it the structure of a one-dimensional
algebraic  group, analogous to $\mathbb{G}_a$ and $\Gm$. By analogy with $\frak h/W =
(\mathbb{G}_a\otimes \Lambda)/W$ and $H/W = (\Gm\otimes \Lambda)/W$, one can also
consider $(E\otimes \Lambda)/W =\mathcal{M}_G(E)$. By a theorem of Looijenga,
$\mathcal{M}_G(E)$ is the coarse moduli space of semistable $G$-bundles over $E$, and
has been extensively studied in \cite{FM, FMII} via the so-called parabolic
construction.  The parabolic construction extends to the case of a nodal or cuspidal
curve and defines a compactification of the adjoint quotient of $G$ or $\frak g$.
Over the adjoint quotient itself, the parabolic construction defines a universal
bundle,  and this is the method used in \cite{FMIII} to define a section of the
adjoint quotient morphism. 

For a smooth elliptic curve
$E$, one can try to use the methods of this paper to investigate different universal
bundles
$\Xi$ over $\mathcal{M}_G(E) \times E$. However, there are some important differences.
The bundle $\Xi$ might only exist as a conformal bundle over $\mathcal{M}_G(E) \times
E$, i.e.\ as a bundle over some quotient of $G$ by a subgroup of the center. If $G$ is
not of type $A_n$ or $C_n$, then $\mathcal{M}_G(E)$ is singular, and the bundle $\Xi$
will not exist, even locally or in some conformal form of $G$, over the singular
points. Finally, even when $\mathcal{M}_G(E)$ is smooth, the bundle $\Xi$ will not be
unique. However, the methods of this paper, in a relative setting, can be used in many
cases to study families of regular semistable bundles over $E$ or over an elliptic
fibration.

The general method for proving our main theorem is as follows. We work entirely in
the context of $G$-bundles. Let $E$ be a Weierstrass cubic, i.e.\ a reduced irreducible
plane cubic curve with a marked smooth  point $p_0$. 
Define
$\widetilde{\mathcal{M}} = E_{\rm reg}\otimes
\Lambda$ and
$\mathcal{M} = \widetilde{\mathcal{M}}/W$. Given a universal bundle
$\Xi\to \mathcal{M}\times E$, there is a corresponding automorphism sheaf
$\mathcal{A}$ over $\mathcal{M}$, defined by  $\mathcal{A}(U) = \Aut (\Xi|U\times E)$.
One can show that, up to isomorphism, $\mathcal{A}$ is independent of the choice of
$\Xi$.  A standard argument shows that the set of all universal bundles is  a principal
homogeneous space over $H^1(\mathcal{M}; \mathcal{A})$. Thus we must show that
$H^1(\mathcal{M}; \mathcal{A})$ is trivial. This is done by analyzing the sheaf
$\mathcal{A}$ and relating it to various finite covers of $\mathcal{M}$ associated to
$\Xi$ as follows. Let $\rho\colon G \to GL(V)$ be a representation of $G$, which in
practice will always be minuscule or quasi-minuscule (to be defined in Section 1). We
can then form the associated vector bundle
$\mathcal{V}= \Xi\times_GV$, which is a family of semistable vector bundles over
$\mathcal{M}\times E$. The Fourier-Mukai correspondence functorially associates to
$\mathcal{V}$ a finite cover $T$ of $\mathcal{M}$, not necessarily reduced, and a sheaf
of $\scrO_T$-modules $\mathcal{Q}$, finite and flat over $\mathcal{M}$, which
determines $\mathcal{V}$. The scheme $T$ is called the \textsl{spectral cover} of
$\mathcal{M}$ associated to $\mathcal{V}$, and we are able to describe $\mathcal{A}$
in terms of the sheaf of invertible functions on $T$ and various correspondences
associated to $T$. Once we have this description, we are able to prove the desired
vanishing for $H^1(\mathcal{M}; \mathcal{A})$.

Automorphism sheaves and spectral covers have been studied by many authors, including
Donagi \cite{Donagi}, Donagi-Gaitsgory \cite{DG}, and Kanev \cite{Kanev}. Knop has
shown (cf.\ \cite{Knop}) that, if $\pi\colon \widetilde{\mathcal{M}} \to \mathcal{M}$
is the projection and $\underline{H}$ is the sheaf over $\widetilde{\mathcal{M}}$
associated to the constant group scheme $\widetilde{\mathcal{M}} \times H$, then
$\mathcal{A} = (\pi_*\underline{H})^W = (\pi_*\mathbb{G}_{m,
\widetilde{\mathcal{M}}}\otimes \Lambda)^W$. It is likely that this description could
also lead to a proof of the main theorem.

\subsection*{Notation and conventions}

All schemes in this paper will be of finite type over $\Cee$, all morphisms are
morphisms of schemes, all sheaves will be in the \'etale topology and cohomology will
always refer to
\'etale cohomology.  Minor modifications of our methods will prove the analogous results
for the uniqueness  of holomorphic sections of the adjoint quotient. In general, we
will consider sheaves of abelian groups, usually written multiplicatively. The sheaf
$\mathbb{G}_{m,T}$ refers to the sheaf (in the
\'etale topology) of sections of the constant group scheme $T\times \Gm$. The
corresponding sheaf in the Zariski topology will be demoted by $\scrO_T^*$. Given a
coherent sheaf
$\mathcal{F}$ over a scheme $X$, we will also use $\mathcal{F}$ to denote the
associated sheaf in the
\'etale topology. We will usually identify a vector bundle with its locally free
sheaf of sections. Finally,
$E$ always denotes a Weierstrass cubic, necessarily either cuspidal, nodal,  or smooth,
with a marked smooth point $p_0$.

\section{Preliminaries}

\subsection{Preliminaries on sheaves}

\begin{defn}  Let
$X$ be a reduced irreducible scheme and let $\cS$ be a sheaf of sets or abelian groups
on
$\cS$. 
\begin{enumerate}
\item[\rm (i)] The sheaf $\cS$ is \textsl{torsion free} if, for every \'etale
open set $f\colon U\to X$,   all sections $s, s'\in \Gamma(U,\cS)$ and every Zariski
open dense subset $U'$ of $U$, if $s|U'=s'|U'$, then $s=s'$. Equivalently, for every
nonempty Zariski open subset
$V$ of $X$, if $j\colon V \to X$ is the inclusion, then the natural morphism
$\cS\to j_*j^*\cS$ is an injection. Since $\cS$ is a sheaf in the \'etale topology, for
every
\'etale open set $U\to X$,  if $U'\to U$ is an \'etale open set with dense image and 
$s, s'\in \Gamma(U,\cS)$ are such that their  restrictions to $U'$ are equal, then
$s=s'$.   Note that every subsheaf of a torsion free sheaf is torsion free.
\item[\rm (ii)] The sheaf $\cS$ has the \textsl{Hartogs property} if for every
\'etale open set  $U\to X$ and every Zariski open subset
$U'\subseteq U$ whose complement is of codimension at least $2$, 
every section $s\in \Gamma(U',\cS)$ extends to a section $\hat s\in \Gamma(U,\cS)$. If
in addition $\cS$ is torsion free,  then there is a unique such $\hat s$. Equivalently,
$\cS$ has the Hartogs property if for every Zariski open subset $V$ of $X$ such that
$X-V$ has codimension at least two, if $j\colon V \to X$ is the inclusion, then
$\cS\to j_*j^*\cS$ is surjective, and is an isomorphism if in addition
$\cS$ is torsion free.
\item[\rm (iii)] If $\cS$ is a torsion free sheaf of sets over $X$ and
$\cS'\subseteq
\cS$ is a subsheaf, then   $\cS'$ is a \textsl{closed subsheaf} of $\cS$ if for every
\'etale open set  $f\colon U\to X$, every section $s\in \Gamma(U,\cS)$ and every \'etale
open dense subset $U'$ of $U$, if $s|U'$ lies in  $\Gamma(U',\cS')$, then $s\in
\Gamma(U,\cS')$.  A standard argument shows that it is enough to check this property
for every Zariski open and dense subset $U'$ of $U$. If
$\cS$ is a sheaf of abelian groups and
$\cS'$ is a subsheaf of abelian groups, then $\cS'$ is a   closed subsheaf of
$\cS$ if and only if the quotient sheaf $\cS/\cS'$ is torsion free.
\end{enumerate}
\end{defn}

For example, if $V'$ is a subbundle of a vector bundle $V$, then the sheaf of sections
of  $V'$ in the \'etale topology is a closed subsheaf of the sheaf of sections of $V$.

The following lemma is then clear:

\begin{lemma}\label{Hartogssub} \begin{enumerate}
\item[\rm (i)] Suppose that $\cS$ is torsion free and
has the Hartogs property. Let
$\cS'$ be a closed subsheaf of $\cS$. Then $\cS'$ is torsion free and has the Hartogs
property as well.
\item[\rm (ii)] Suppose that $\cS'$ is a subsheaf of a torsion free sheaf $\cS$, and
that
$\mathcal{C}$ is a closed subsheaf of $\cS$. Then $\mathcal{C}\cap\cS'$ is a closed
subsheaf of
$\cS'$.
\item[\rm (iii)] If $\cS$ and $\cS'$ are torsion free sheaves of abelian groups and
$c\colon \cS \to \cS'$ is a homomorphism, then $\Ker c$ is a closed subsheaf of $\cS$.
\end{enumerate}
\qed
\end{lemma}

\subsection{A geometric lemma}

\begin{lemma}\label{geomlemma}  Let $X$ be a scheme and $T$ a torus such that there
exists a morphism
$\pi\colon X
\to T$ and a finite \'etale covering $\widetilde T\to T$ such that the fiber
product
$X\times _{T}\widetilde T$ is isomorphic to $\Aff^k\times \widetilde T$ as schemes over
$\widetilde T$. Then
\begin{enumerate}
\item[\rm (i)] Every morphism $f\colon X \to \Gm$ is induced by a morphism $T\to \Gm$,
and hence is given up to translation by an element of $\Gm$ by a character of $T$. In
particular, every morphism from an affine space to $\Gm$ is constant.
\item[\rm (ii)] $\Pic X$ is trivial.
\end{enumerate}
\end{lemma}
\begin{proof} The proof of (i) is left to the reader. To see (ii), let $L$ be a line
bundle on $X$ and let $d$ be the degree of the morphism $\widetilde T \to T$.  Then the
pullback $\widetilde L$ of
$L$ to $X\times _{T}\widetilde T$ is trivial, as is $\operatorname{Norm}(\widetilde L)
= L^d$. Thus every element in $\Pic X$ is $d$-torsion. Let $(\Pic X)_d$ be the group
of $d$-torsion line bundles of $X$. The Kummer exact sequence gives an exact sequence
$$H^0(X; \Gm) \xrightarrow{z\mapsto z^d} H^0(X; \Gm) \to H^1(X; \boldsymbol{\mu}_d)\to
(\Pic X)_d
\to 0.$$ Since $X$ is a fiber bundle over $T$ with fibers isomorphic to $\Aff^k$, the
comparison of \'etale and classical cohomology implies that
$H^1(X;
\boldsymbol{\mu}_d)
\cong (\Zee/d\Zee)^r$, where
$r=\dim T$. By Part (i), 
$H^0(X; \Gm)\cong \Cee^*\times \mathsf{X}(T)$, where $\mathsf{X}(T)$, the character
group of
$T$, is a free abelian group of rank $r$. Since the homomorphism $H^0(X; \Gm) \to
H^0(X; \Gm)$ is raising to the power $d$, it follows that the cokernel has order $d^r$.
Hence the cokernel maps isomorphically onto $H^1(X; \boldsymbol{\mu}_d)$ and $(\Pic
X)_d$ is trivial. Thus $\Pic X$ is trivial as well.
\end{proof}

\subsection{Minuscule and quasi-minuscule weights}

We begin by recalling some standard definitions from representation theory. The
following three sets of objects may be identified:
homomorphisms $\Lambda \to \Zee$,  linear maps $\frak h\to \Cee$ which are integral
on $\Lambda$, and  characters $H\to \Cee^*$. We will call a homomorphism 
$\mu\colon \Lambda
\to \Zee$  a
\textsl{weight} and will also $\mu$ to denote the obvious extension of $\mu$ to a linear
map
$\frak h \to \Cee$ and occasionally the corresponding character $H\to \Cee^*$.

\begin{defn} Let $\rho\colon G\to GL(V)$ be an irreducible representation with highest
weight $\varpi$. Then
$\rho$ or $\varpi$ is \textsl{minuscule} if for every pair of weights $\mu_1, \mu_2$ of
$\rho$, there exists a $w\in W$ such that $w\mu_1=\mu_2$. 
The representation $\rho$ is \textsl{quasi-minuscule} if for every pair of
non-zero weights
$\mu_1, \mu_2$ of
$\rho$, there exists a $w\in W$ such that $w\mu_1=\mu_2$. Since highest weights have
multiplicity one, if
$\rho$ is quasi-minuscule, every non-zero weight of $\rho$ has multiplicity one.
Finally,
$\rho$ is \textsl{properly quasi-minuscule} if it is quasi-minuscule but not minuscule.
\end{defn}

The following is a standard result, cf.\ \cite{Bour}:

\begin{lemma} With notation as above,
\begin{enumerate}
\item[\rm (i)] The minuscule representations are exactly the representations whose
highest weight is a fundamental weight $\varpi_\alpha$ for a simple root $\alpha$ such
that, if $\widetilde \alpha =\sum _{\beta \in \Delta}h_\beta\beta$ is the highest root,
then
$h_\alpha = 1$. Thus, the number of isomorphism classes of  minuscule representations is
$\#Z(G) -1$, where
$Z(G)$ is the center of $G$. 
\item[\rm (ii)] For each $G$, there is exactly one properly quasi-minuscule
representation $\rho$ up to isomorphism. Its highest weight is the highest short root of
$G$. Hence, if
$G$ is simply laced, $\rho$ is the adjoint representation.
\qed
\end{enumerate}
\end{lemma}

\begin{remark}\label{qminforms} Let $\rho\colon G \to GL(V)$ be a
properly quasi-minuscule representation. Let $V=\bigoplus V_\mu$ be the
decomposition of $V$ into weight spaces, where the $\mu$ are either $0$ or the short
roots of $G$ and, for
$\mu\neq 0$, $\dim V_\mu =1$.  Then  there exists a non-degenerate $G$-invariant
symmetric bilinear form
$B\colon V\otimes V \to \Cee$. If $\mu_1$ and $\mu_2$ are nonzero, then $B(V_{\mu_1},
V_{\mu_2}) \neq 0$ if and only if $\mu_1+\mu_2 =0$.
There also exists a   $G$-invariant $3$-tensor
$C\colon V\otimes V\otimes V\to \Cee$, with the property that, if $\mu_1, \mu_2, \mu_3$
are nonzero, then
$C(V_{\mu_1}, V_{\mu_2}, V_{\mu_3}) \neq 0$ if and only if $\mu_1+\mu_2 +\mu_3=0$.
For example, if $G$ is simply laced, then $\rho$ is the adjoint representation. In
this case, we let $B(X,Y)$ be the Killing form and let $C(X,Y,Z)$ be the Cartan form
$B([X,Y],Z)$. It is easy to see that $B$ and $C$ have the desired properties. The
remaining groups are done on a case-by-case basis. For $G=Spin (2n+1)$, the
representation $\rho$ is the standard representation $G \to SO(2n+1) \subseteq GL(2n+1,
\Cee)$. In this case, $B$ is the standard form and $C=0$. (Note that two short roots of
$G$ never sum to a short root.) If $G=Sp(2n, \Cee)$, let $V_0$ be the standard
$2n$-dimensional representation of $G$ with the invariant symplectic form $\langle
\cdot, \cdot\rangle$. Then $V$ is a direct summand of $\bigwedge^2V_0$ complementary to
the line defined by $\langle
\cdot, \cdot\rangle$. The form $B$ is then the restriction to $V$ of the form on
$\bigwedge^2V_0$ induced by
$$B(v_1\wedge v_2, w_1\wedge w_2) = \langle v_1, w_1\rangle\langle v_2, w_2\rangle
+\langle v_1, w_2\rangle\langle v_2, w_1\rangle.$$
Viewing $C$ as a $G$-invariant homomorphism $\varphi\colon V\otimes V \to V^* \cong V$
via the form
$B$, the tensor $C$ corresponds to
$$\varphi(v_1\wedge v_2, w_1\wedge w_2) = \langle v_2, w_1\rangle v_1\wedge w_2-\langle
v_1, w_1\rangle  v_2\wedge w_2 -\langle v_2, w_2\rangle v_1\wedge w_1 + \langle v_1,
w_2\rangle v_2\wedge  w_1.$$
The definitions of $B$ and $C$ for the groups of type $G_2$ and $F_4$ are well-known
\cite{Adams}. 

In all cases, $\rho(G)$, which is the adjoint form of $G$, is exactly the stabilizer
in $GL(V)$ of the forms $B$ and $C$. 

While we will not use this remark, it motivates the constructions of Section 5.
\end{remark}

\section{Automorphism sheaves and $G$-bundles}

\subsection{Definition of automorphism sheaves}

Let $X$ and $Y$ be   schemes, with $Y$ projective, and let
$\Xi\to X\times Y$ be a principal
$G$-bundle over
$X\times Y$. We define a sheaf $\pi_1{}_*\Aut \Xi$ over $X$ as follows: For each
\'etale open set
$f\colon U\to X$, let $\Xi_U = (f\times \Id)^*\Xi$ and let
$\pi_1{}_*\Aut \Xi(U) = \Aut \Xi_U$. Thus $\pi_1{}_*\Aut \Xi(U)$ is
the group of sections of $\Ad \Xi_U=\Xi_U\times _GG$, where $G$
acts on itself by conjugation. If
$\mathcal{V}$ is a vector bundle over $X\times Y$, the sheaf $\pi_1{}_*\Aut \mathcal{V}$
is defined similarly.

There is also the related functor $\mathbf{A}$ from schemes over $X$ to groups defined
similarly  for each $f\colon X'\to X$ by $\mathbf{A}(X') =\Aut (f\times \Id)^*\Xi$. In
this case, one can show the following:
 
 \begin{proposition}\label{conj} Let $E$ be a singular Weierstrass cubic. Suppose  that
$\Xi\to X\times E$ is a principal
$G$-bundle such that, for every $x\in X$, $\Xi|\{x\}\times
E$ pulls back to the trivial bundle on $\widetilde E$ and is regular. Then $\mathbf{A}$
is representable by a commutative affine  group scheme which is smooth over $X$. \qed
\end{proposition}

Since we will not use Proposition~\ref{conj} in this paper, we omit the proof.

We now discuss some general properties of $\pi_1{}_*\Aut \Xi$ and $\pi_1{}_*\Aut
\mathcal{V}$. First we have the following:

\begin{lemma} Let $X$ be normal and let $\mathcal{V}$ be a locally free sheaf over
$X\times E$. Then $\pi_1{}_*\mathcal{V}$ has the Hartogs property.
\end{lemma}
\begin{proof} It clearly suffices to show that $\scrO_{X\times E}$ has the Hartogs
property. This is an easy consequence of the fact that $X\times \widetilde{E}$ is
normal and that, locally,  any section of $\scrO_{X\times \widetilde{E}}$ which lies in
$\scrO_{X\times E}$ in the complement of a closed subset of codimension at least two
lies in $\scrO_{X\times E}$.
\end{proof}

\begin{lemma}\label{AutHartogs} Let $\mathcal{V}$ be a vector bundle over $X\times E$,
where
$X$ is normal. Then $\pi_1{}_*\Aut \mathcal{V}$ is torsion free and has the Hartogs
property.
\end{lemma}
\begin{proof} The sheaf $\pi_1{}_*\Aut \mathcal{V}$ is the subsheaf of invertible
elements in the sheaf of algebras  $\pi_1{}_*Hom( \mathcal{V},
\mathcal{V})$. Since
$\mathcal{V}$ and hence $Hom( \mathcal{V},
\mathcal{V})$ are vector bundles and $X$ is normal, $\pi_1{}_*Hom( \mathcal{V},
\mathcal{V})$ is torsion free and has the Hartogs property. Hence $\pi_1{}_*\Aut
\mathcal{V}$ is torsion free. Suppose that $U\to X$ is an
\'etale open set and $U'=U-Y$ where $Y$ has codimension at least two in $U$. If
$s\in \Gamma(U',\pi_1{}_*\Aut \Xi)$, then $s^{-1} \in \Gamma(U',\pi_1{}_*\Aut \Xi)$
as well, and both $s$ and $s^{-1}$ extend to elements of  $\Gamma(U, \pi_1{}_*Hom(
\mathcal{V},
\mathcal{V}))$ whose product is $\Id$ on $U'$ and hence on $U$. Thus the unique
extension of $s$ lies in 
$\Gamma(U, \pi_1{}_*\Aut \mathcal{V})$, so that $\pi_1{}_*\Aut \mathcal{V}$ has the
Hartogs property.
\end{proof}

\begin{lemma}\label{general} Let $\Xi\to X\times E$ be a
principal $G$-bundle over $X\times E$, and suppose that $X$ is normal.
\begin{enumerate}
\item[\rm (i)] Let $\rho\colon G\to \Cee^N$ be a faithful representation and let
$\mathcal{V} =\Xi\times _G\Cee^N$. Then $\pi_1{}_*\Aut \Xi$ is a closed subsheaf of
$\pi_1{}_*\Aut \mathcal{V}$.
\item[\rm (ii)] $\pi_1{}_*\Aut \Xi$ is torsion free and has the Hartogs property.
\end{enumerate}
\end{lemma}
\begin{proof} To see (i), let $f\colon U\to X$ be an  \'etale open set  and let
$\Xi_U=(f\times \Id)^*\Xi$. The scheme $\Ad \Xi_U=\Xi_U\times _GG$ is a closed
subscheme of
$\Xi_U\times_GGL(N,\Cee)$. Thus, if $U'\to U$ is a \'etale open and dense subset of
$U$ and
$s$ is a section of $\Xi_U\times_GGL(N,\Cee)$   whose restriction to $U'\times E$ lies
in $\Ad \Xi_{U'}$, then   $s$ lies in $\Ad\Xi_U$. This proves that 
$\pi_1{}_*\Aut \Xi$ is a closed subsheaf of
$\pi_1{}_*\Aut \mathcal{V}$. Since $\pi_1{}_*\Aut \mathcal{V}$ is torsion free,
$\pi_1{}_*\Aut \Xi$ is torsion free as well, and by Lemma~\ref{Hartogssub} and
Lemma~\ref{AutHartogs},
$\pi_1{}_*\Aut \Xi$ has the Hartogs property.
\end{proof}

\subsection{Statement of the main theorem}

We now specialize to the case of a singular Weierstrass cubic $E$. Recall the notation
of the introduction: let $\widetilde{\mathcal{M}} = E_{\rm reg}\otimes \Lambda$, so
that $\widetilde{\mathcal{M}}  = \frak h$ if $E$ is cuspidal and
$\widetilde{\mathcal{M}} = H$ if $E$ is nodal, and let $\mathcal{M}
=\widetilde{\mathcal{M}}/W$ be the adjoint quotient. 

\begin{defn} Fix once and for all a section of
the adjoint quotient morphism, for example the Kostant or Steinberg section, and let
$\Xi\to
\mathcal{M}\times E$ be the corresponding
$G$-bundle.  Let $\mathcal{A} =\pi_1{}_*\Aut \Xi$, where $\Xi\to \mathcal{M}\times
E$ is a universal $G$-bundle as above.
\end{defn}

Since the centralizer of a regular element is abelian, $\mathcal{A}$ is a sheaf of
abelian groups over $\mathcal{M}$. It is easy to see that two different choices of a
section lead to isomorphic automorphism sheaves. Note further that there is a natural
inclusion of the center
$Z(G)$ of
$G$, viewed as a constant sheaf on $\mathcal{M}$, into $\mathcal{A}$.

If
$\sigma'$ is another  section of the adjoint
quotient morphism, let
$\Xi'$ be the corresponding bundle. Then, for all
$x\in
\mathcal{M}$, there exists a
$g\in G$ such that $\ad(g)\sigma(x) =\sigma'(x)$, and any two such differ by right
multiplication by
$Z_G(\sigma(x))$. Since $\sigma$ is a section of the adjoint quotient morphism, the
morphism $G\times 
\frak h/W \to \frak g_{\rm reg}$ defined by $(g,x)\mapsto \ad g(\sigma(x))$ is smooth.
Thus there exist local sections of this morphism in the \'etale topology. Similarly in
the case of the adjoint quotient of $G$,   the
morphism $G\times 
H/W \to G_{\rm reg}$ defined by $(g,x) \mapsto \Ad g(\sigma (x))$ is smooth. A formal
argument then shows: 

\begin{proposition}\label{formal} Let $\sigma$ and $\sigma'$ be two sections of the
adjoint quotient morphism $\frak g_{\rm reg} \to \frak h/W=\mathcal{M}$ or $G_{\rm
reg}  \to H/W =\mathcal{M}$. For every
$x\in \mathcal{M}$, there exists an \'etale neighborhood $U$ of $x$ and a function
$\gamma\colon U\to G$ such that $\sigma'|U =\ad(\gamma)(\sigma|U)$, in the case of
$\frak g$, or $\sigma'|U =\Ad(\gamma)(\sigma|U)$ in the case of $G$. Moreover, any two
choices of
$\gamma$ differ by right multiplication by an element of
$\mathcal{A}(U)$.
\qed
\end{proposition}

In the standard way, choosing a sufficiently small open cover $\{U_i\}$ of $\mathcal{M}$
and writing
$\sigma'|U_i = \ad (\gamma_i)(\sigma|U_i)$ with $\gamma _j|U_i\times _{\mathcal{M}}U_j
= \gamma_ia_{ij}|U_i\times _{\mathcal{M}}U_j$ defines a $1$-cocycle $\{a_{ij}\}$ and
hence an element of  
$H^1(\mathcal{M};
\mathcal{A})$ which is independent of the choices. This element is the identity if and
only if the bundles
$\Xi$ and
$\Xi'$ are isomorphic, if and only if there exists a morphism $g\colon \mathcal{M}\to
G$ such that
$\sigma'=\ad g(\sigma)$, in which case we say that the sections $\sigma$ and $\sigma'$
are \textsl{conjugate}. This explains the relevance of the second half of the following
theorem:

\begin{theorem}\label{main} Let $\mathcal{M}$ be the adjoint quotient of either $\frak
g$ or
$G$, and let $\mathcal{A}$ denote the automorphism sheaf corresponding to the section
$\sigma$ with image contained in the set of regular elements. Then 
$$H^0(\mathcal{M}; \mathcal{A}) = \begin{cases} Z(G), &\text{for the case of $\frak
g$;}\\ Z(G)\times \Zee , &\text{for the case of $G$.}
\end{cases}$$
Moreover, in both cases, $H^1(\mathcal{M}; \mathcal{A}) =\{1\}$.
\end{theorem}

\begin{corollary} Every two sections of the adjoint quotient morphism are conjugate by
a morphism from the adjoint quotient to the set of regular elements of $\frak g$ or
$G$. \qed
\end{corollary}

\section{Quasi-minuscule representations and spectral covers}

\subsection{Vector bundles attached to   quasi-minuscule representations}

 We fix the following notation:
 $\rho\colon G\to \Aut (V)$
is a  quasi-minuscule representation, with highest weight $\varpi$ (with
respect to the given fixed set of simple roots). Let
$W_0\subseteq W$ be the stabilizer of $\varpi$. The
group $W_0$ is the group generated by reflections in the roots
annihilated by $\varpi$. 

With $\mathcal{M}$ as before, and $\Xi$ the universal bundle
over $\mathcal{M}\times E$ corresponding to the choice of the section $\sigma$, 
let ${\mathcal V}\to {\mathcal M}\times E$ be the
vector bundle induced by the representation $\rho$. By \cite[Theorem 2.4.3]{Spectral},
we have:

\begin{proposition} The sheaf $\pi_1{}_*{\mathcal V}$ is a locally free
sheaf over $\mathcal{M}$ whose rank is the multiplicity of the trivial
weight in $\rho$. Moreover, the inclusion $\pi_1{}^*\pi_1{}_*{\mathcal V} \to
\mathcal{V}$ is a subbundle, so that the quotient  $\overline
{\mathcal V}$ is a vector bundle, and there is an exact sequence of vector bundles over
${\mathcal M}\times E$:
$$0\to \pi_1{}^*\pi_1{}_*{\mathcal V}\to
{\mathcal V}\to \overline{\mathcal V}\to 0.\qed$$
\end{proposition}

We shall need the following facts about $\Xi$ and $\mathcal{V}$:

\begin{lemma}\label{Gsplit}
There is  an \'etale open and dense subset  $U\to {\mathcal M}$
such that
$\Xi_U$ has a reduction of structure to an $H$-bundle $\Xi_{U,H}$. The action of the
constant group scheme $U\times H$ over $U$ as a group scheme of automorphisms of
$\Xi_{U,H}$ defines an isomorphism from the sheaf of sections of the constant group
scheme
$U\times H\to U$ to
${\mathcal A}|U$. 
\end{lemma}

\begin{proof}
Let $U$ be the nonempty open subset of $\widetilde {\mathcal M}$
where no root is the identity. Thus $U$ is either $\frak h_{\rm reg} =\frak h\cap \frak
g_{\rm reg}$ or $H_{\rm reg} = H\cap G_{\rm reg}$. We just write out the Lie algebra
case. The induced morphism
$U\to {\mathcal M}$ is \'etale. For a scheme $B$, an $H$-bundle over $B\times E$
which is trivialized over $B\times \widetilde{E}$ is identified with a morphism $B\to
\frak h$. Let  
$\Xi_{H,U}$ be the
$H$-bundle over
$U\times E$ induced by applying this remark to  the inclusion $U=\frak h_{\rm
reg}\to \frak h$ and let $\Xi'_U = \Xi_{H,U}\times_HG$. It is clear that $\Xi_U$ and
$\Xi_U'$ have isomorphic restriction to every slice. A formal argument as in the proof
of Proposition~\ref{formal} shows that, possibly after replacing $U$ by  an \'etale
open dense subset, the bundles $\Xi_U$ and
$\Xi_U'$ become isomorphic. This gives the desired reduction of structure.   Viewing
$\Xi_U \cong \Xi_U'$ as the
$G$-bundle associated to a morphism $\sigma'\colon U \to \frak h_{\rm reg} \subseteq
\frak g$, the sections of $\mathcal{A}$ over an \'etale open set $V$ of $U$ are given by
morphisms
$f\colon V \to G$ such that $\ad (f)(\sigma'|V) = \sigma'|V$.  Clearly, since the
image of $\sigma'$ is contained in $\frak h_{\rm reg}$, this is the same as the group of
morphisms from
$V$ to
$H$.
\end{proof}

\begin{corollary}\label{split}
Let $U \to \mathcal{M}$ be the \'etale open subset of Lemma~\ref{Gsplit}. Then 
${\mathcal V}|U\times E$ is isomorphic to
$\bigoplus_{\mu}\mathcal{L}_\mu^{\oplus m_\mu}$, where the $\mu$ are the weights of
$\rho$,
$m_\mu$ is the multiplicity of $\mu$ in
$\rho$, and the $\mathcal{L}_\mu$ are line bundles on $U\times E$ with the property that
$\mathcal{L}_{\mu}\cong \mathcal{L}_{\mu'}$ if and only if $\mu=\mu'$. Under the
isomorphism of $\mathcal{A}|U$ with the sheaf of sections of the constant group scheme
$U\times H$, a local section $s$ of $U\times H$ acts on $\mathcal{L}_\mu^{\oplus m_\mu}$
via multiplication by  the function $\mu(s)$.
\qed
\end{corollary}

As before, let $\mathcal A$ be the automorphism sheaf of $\Xi$.   The representation
$\rho$ induces a homomorphism ${\mathcal A}\to \pi_1{}_*{\Aut}({\mathcal V})$. The
kernel of this homomorphism is the subsheaf of constant automorphisms with values
in   ${\Ker}(\rho)\subseteq Z(G)$. Let $\overline{\mathcal
A}_\rho$ denote the quotient sheaf ${\mathcal A}/\Ker(\rho)$.

\begin{lemma}\label{AvsAbar} In the above notation, there is an exact sequence
$$\{1\}\to \Ker(\rho) \to  H^0(\mathcal{M}; \mathcal A) \to H^0(\mathcal{M}; \overline
{\mathcal A}_\rho) \to \{1\}$$ and the natural map $H^1(\mathcal{M};
\mathcal A) \to  H^1(\mathcal{M}; \overline {\mathcal A}_\rho)$ is an isomorphism.
\end{lemma}
\begin{proof} This is an immediate consequence of the fact that, as $\mathcal{M}$ is
isomorphic to $\Aff^r$, $H^0(\mathcal{M};\Ker(\rho))\cong \Ker(\rho)$ and
$H^i(\mathcal{M};\Ker(\rho))$ is trivial for $i=1,2$.
\end{proof}

\begin{proposition}\label{closure} The sheaves
$\mathcal A$ and $\overline{\mathcal A}_\rho$ are torsion free and have the Hartogs
property. The sheaf $\overline{\mathcal A}_\rho$ is a closed subsheaf of
$\pi_1{}_* \Aut({\mathcal V})$. Finally, every local
section of $\overline{\mathcal A}_\rho$ induces the identity on
$R^0\pi_1{}_*{\mathcal V}$.
\end{proposition}

\begin{proof} The first two statements follow from Lemma~\ref{general}. It suffices to
prove the final statement on an \'etale open dense of $\mathcal{M}$. By
Corollary~\ref{split}, it is immediate that $\overline{\mathcal A}_\rho$ has this
property on the \'etale open set 
$U\to{\mathcal M}$   of
Lemma~\ref{Gsplit}. 
\end{proof}

\subsection{The spectral cover}

We begin by recalling the general formalism of the Fourier-Mukai correspondence in the
spcecial case we shall need, for which a general reference is \cite[\S1]{Spectral}. Let
$B$ be a scheme  and let
$\mathcal{V}$ be a vector bundle of rank $n$ over
$B\times E$ whose restriction to every slice $\{x\}\times E$ has trivial pullback to the
normalization of $E$. Then $\mathcal{V}$ determines functorially  a relative Cartier
divisor
$T_\mathcal{V}\subseteq B\times E_{\rm reg}$, finite over
$B$, and a sheaf $\mathcal{Q}_\mathcal{V}$ of $\scrO_{T_\mathcal{V}}$-modules, finite
and flat over
$B$ of degree
$n$, such that $\mathcal{V} =q_*(p_1^*(\mathcal{Q}_\mathcal{V}\otimes
\pi_2^*\scrO_E(p_0))\otimes r^*\mathcal{P})$, where $p_1\colon B\times E\times E$ is
projection onto the first two factors, $q\colon B\times E\times E$ is projection onto
the first and third factors, $r\colon B\times E\times E\to E\times E$ is projection
onto the last two factors, and $\pi_2\colon B\times E\to E$ is projection onto the
second factor. The divisor $T_\mathcal{V}$ is additive under exact sequences, and,
viewing $\mathcal{Q}_\mathcal{V}$ as a sheaf on $B\times E$, given an exact sequence 
$0\to \mathcal{V}' \to \mathcal{V} \to \mathcal{V}'' \to 0$ of vector bundles as above,
there is an exact sequence of the corresponding sheaves:
$$0\to \mathcal{Q}_{\mathcal{V}'} \to \mathcal{Q}_{\mathcal{V}} \to
\mathcal{Q}_{\mathcal{V}''} \to 0.$$
Moreover the constructions of $T_\mathcal{V}$ and $\mathcal{Q}_\mathcal{V}$ are
compatible with base change. We somewhat loosely refer to
$T_\mathcal{V}$ as the
\textsl{spectral cover} of
$B$ corresponding to $\mathcal{V}$. More generally, we have the following:

\begin{defn} Let $\nu\colon T \to B$ be a finite
morphism of schemes.  Suppose that
$\mathcal{Q}$ is a coherent sheaf on $T$, flat over $B$, and that $\mathcal{L}$ is a
line bundle on $T\times E$. The \textsl{spectral cover construction} (applied to
$\mathcal{Q}$ and $\mathcal{L}$) is the vector bundle 
$$\mathcal{V}= (\nu\times \Id)_*(\pi_1^*\mathcal{Q}\otimes \mathcal{L}).$$
\end{defn}

In the above situation, the sheaf $\nu_*\scrO_T$ acts via endomorphisms on 
$\mathcal{V}$, and hence $\nu_*\scrO_T^*$ acts via automorphisms of 
$\mathcal{V}$. Passing from a sheaf in the Zariski topology to one in  the \'etale
topology defines a homomorphism
$\nu_*\mathbb{G}_{m,T}\to \pi_1{}_*\Aut \mathcal{V}$.

\begin{lemma}\label{Gmclosed} In the above notation, suppose that $B$ is
reduced irreducible and that $\mathcal{Q}$ is a locally free
$\scrO_T$-module, or more generally that $T$ is flat over $B$ and that, for all $x\in
B$,  
$H^0(\nu^{-1}(x); \mathcal{Q}|\nu^{-1}(x))$ is a faithful
$H^0(\nu^{-1}(x);\scrO_T|\nu^{-1}(x))$-module. Then
$\nu_*\mathbb{G}_{m,T}$ is a closed subsheaf of
$\pi_1{}_*\Aut \mathcal{V}$.
\end{lemma}
\begin{proof} The fiber of $\mathcal{V}$ over a point $(x,e)$ of $B\times E$ is
identified up to a homothety with $H^0(\nu^{-1}(x); \mathcal{Q}|\nu^{-1}(x))$ and the
induced action of $\nu_*\scrO_T$ on this fiber is via the action of
$H^0(\nu^{-1}(x);\scrO_T|\nu^{-1}(x))$. Since $T$ is flat over $B$, $\nu_*\scrO_T$
is a vector bundle. The hypotheses imply that, for each $x\in B$, $(\nu_*\scrO_T)_x \to
\pi_1{}_*Hom(\mathcal{V}, \mathcal{V})_x$ is injective. Thus
$\nu_*\scrO_T$ is a vector subbundle, and hence a closed
subsheaf,  of
$\pi_1{}_*Hom(\mathcal{V},\mathcal{V})$.  By (i) of Lemma~\ref{Hartogssub}, 
as sheaves of sets in the \'etale topology, $\pi_1{}_*\Aut\mathcal{V}\cap \nu_*\scrO_T$
is a closed subsheaf of
$\pi_1{}_*\Aut\mathcal{V}$.  It suffices to prove that $\pi_1{}_*\Aut\mathcal{V}\cap
\nu_*\scrO_T= \nu_*\mathbb{G}_{m,T}$. Replacing $B$ by an \'etale open set and $T$ by
the fiber product, it is enough to show that, if $f\in  \Gamma(T,\scrO_T)$ defines an
automorphism of $\mathcal{V}$ and hence by functoriality of $\mathcal{Q}$, then,  for
all
$y\in T$, $f$ does not lie in the maximal ideal $\frak m_y$ of $y$. Since $f$ defines
an automorphism of $\mathcal{Q}$, it defines an automorphism of the stalk of
$\mathcal{Q}$ over $y$, which is a finite module over the local ring of $T$ at $y$ and
is nonzero by our assumptions on $\mathcal{Q}$. It follows by Nakayama's lemma that
$f\notin
\frak m_y$.  Thus
$\pi_1{}_*\Aut
\mathcal{V}\cap \nu_*\scrO_T = \nu_*\mathbb{G}_{m,T}$ is a
closed subsheaf of
$\pi_1{}_*\Aut \mathcal{V}$.
\end{proof}

\subsection{The case of a quasi-minuscule representation}

We return to the case where $\mathcal{V}=\Xi\times_G\Cee^N$ is the vector bundle over
$\mathcal{M}\times E$ defined by the minuscule or quasi-minuscule representation
$\rho\colon G\to GL(N,\Cee)$.  Recall that we have the exact sequence
$$0 \to \pi_1{}^*\pi_1{}_*\mathcal{V}\to \mathcal{V}\to \overline{\mathcal{V}} \to
0.$$
Let
$\hat{T}_0 =
\widetilde{\mathcal{M}}/W_0$ and let $\hat\nu\colon \hat{T}_0\to \mathcal{M}$ be the
natural morphism. The homomorphism
$\varpi\colon
\Lambda \to \Zee$ defines a $W_0$-invariant homomorphism $\widetilde{\mathcal{M}} \to
E_{\rm reg}$ and hence a morphism $f_{\varpi} \colon \hat{T}_0\to E_{\rm reg}$. Let
$T_0\subseteq \mathcal{M}\times E_{\rm reg}$ be the image of the product morphism
$(\hat\nu, f_{\varpi})$. By \cite[Corollary 2.2.3]{Spectral}, the morphism $\hat{T}_0\to
T_0$ is birational and exhibits $\hat{T}_0$ as the normalization of $T_0$. Moreover,
by \cite[Lemma 2.1.1]{Spectral}, we have:

\begin{lemma}\label{hatT1} Let $\Lambda_0 = \Ker \varpi$ and let $\mathcal{M}_0 =
(E_{\rm reg}\otimes
\Lambda_0)/W_0$. There is a finite
\'etale cover
$C\to E_{\rm reg}$, with $C$ necessarily isomorphic to $E_{\rm reg}$, such that
$\hat{T}_0\times _{E_{\rm reg}}C
\cong
\mathcal{M}_0\times C$ as schemes over $C$. \qed
\end{lemma}

Let $T=T_{\mathcal{V}}\subseteq {\mathcal M}\times E_{\rm reg}$ be the   spectral
cover of $\mathcal{V}$, let $\nu\colon T \to \mathcal{M}$ be the morphism induced by
projection onto the first factor, and let
$\mathcal{Q}=\mathcal{Q}_{\mathcal{V}}$ be the corresponding sheaf of
$\scrO_T$-modules. Let $\mathcal{Q}_1$ be the sheaf  corresponding to
$\pi_1{}^*\pi_1{}_*\mathcal{V}$ and let $\overline{\mathcal{Q}}$ be the sheaf
corresponding to $\overline{\mathcal{V}}$. The general formalism of the previous
subsection and \cite[Corollary 2.2.3]{Spectral} imply: 

\begin{proposition} Let $T_0$ be the divisor defined above, let $\overline T_1$ be the
reduced divisor
$\mathcal{M}\times
\{p_0\}$, and let $m_0$ be the multiplicity of the trivial
weight in
$\rho$.  Then the  divisor $T_0$ is the spectral cover of
$\overline{\mathcal{V}}$, and 
$T=T_0\cup T_1\subset {\mathcal M}\times E_{\rm reg}$ where, as a divisor,  
$T_1=m_0\cdot\overline T_1$. Thus $T=T_0$ if and only if $\rho$ is
minuscule.  Generically along $\overline T_1\cap T_0$ the varieties $\overline T_1$ and
$T_0$ are smooth and meet transversely. \qed
\end{proposition}

Let $\overline T$ be the sum  $T_0+ \overline T_1$ considered
as a (reduced) effective Cartier divisor of ${\mathcal M}\times E$. It is the reduced
subscheme of
$T$. By the above, the following is exact on the complement of a codimension two subset
of $\overline{T}$:
$$0\to \scrO_{\overline{T}} \to \scrO_{\overline{T}_1} \oplus \scrO_{T_0} \to
\scrO_{\overline{T}_1\cap T_0} \to 0. $$
In fact, since $\overline{T}$ is Cohen-Macaulay, the above sequence is everywhere
exact, although we will not need this.

\begin{proposition}\label{ovT}
 The sheaf $\mathcal{Q}$ of $\scrO_T$-modules is in
fact a sheaf of $\scrO_{\overline T}$-modules.
\end{proposition}

\begin{proof}
By the formalism of spectral covers, there is an exact sequence
$$0\to  \mathcal{Q}_1\to \mathcal{Q}\to \overline{\mathcal{Q}} \to 0.$$
Since $\pi_1{}^*\pi_1{}_*\mathcal{V}$ is pulled back from a bundle on $\mathcal{M}$,
it restricts to the trivial bundle on every slice and hence $\mathcal{Q}_1$ is
annihilated by $I_{\overline{T}_1}$, the ideal sheaf of $\overline{T}_1$. Since
$\overline{\mathcal{Q}}$ is an $\scrO_{T_0}$-module, it is annihilated by the ideal
sheaf
$I_{T_0}$.  Hence
$\mathcal{Q}$ is annihilated by the product ideal $I_{\overline{T}_1}\cdot I_{T_0}$,
which is the ideal of $\overline T$.
\end{proof}

As for the quotient itself, by \cite[Theorem 2.4.3]{Spectral} and Lemma~\ref{geomlemma},
we have:

\begin{proposition}\label{T1}
The sheaf $\overline{\mathcal{Q}}$ is isomorphic to an invertible
$\scrO_{\hat{T}_0}$-module, and hence is isomorphic, as a sheaf of
$\scrO_{T_0}$-modules, to $\scrO_{\hat{T}_0}$. \qed\end{proposition}

\begin{corollary}\label{AutVbar} $\pi_1{}_*\Aut\overline{\mathcal{V}} = \hat
\nu_*\mathbb{G}_{m, \hat{T}_0}$. In particular, if $\rho$ is minuscule, then
$\mathcal{V} =\overline{\mathcal{V}}$, $\pi_1{}_*\Aut\mathcal{V} = \hat
\nu_*\mathbb{G}_{m, \hat{T}_0}$, and $\overline{\mathcal{A}}_\rho$ is a closed subsheaf
of
$\hat \nu_*\mathbb{G}_{m, \hat{T}_0}$.
\end{corollary}

\begin{proof} By Lemma~\ref{Gmclosed}, $\hat\nu_*\mathbb{G}_{m, \hat{T}_0}$ is a closed
subsheaf of $\pi_1{}_*\Aut\overline{\mathcal{V}}$. Thus, it suffices to prove that they
agree on an \'etale  open and dense  set. By Corollary~\ref{split}, there is an \'etale
open set $U$ such that the pullback of $\mathcal{V}$ to $U$ is of the form
$\bigoplus_\mu \mathcal{L}_\mu^{\oplus m_\mu}$, where the $\mathcal{L}_\mu$ are
distinct line bundles of relative degree zero on
$U\times E$, indexed by the weights of $\rho$, and the $m_\mu$ are the multiplicities
and thus $m_\mu=1$ if $\mu\neq 0$. Hence
$\overline{\mathcal{V}}|U\times E \cong \bigoplus_{\mu\neq 0}\mathcal{L}_{\mu}$. Thus 
$\pi_1{}_*\Aut\overline{\mathcal{V}}|U$ is just a product of copies of
$\mathbb{G}_{m,U}$ indexed by the nonzero weights of $\rho$, and this is clearly equal
to $\hat
\nu_*\mathbb{G}_{m, \hat{T}_0}|U$.
\end{proof}

\subsection{The properly quasi-minuscule case}

In the properly quasi-minuscule case, we will need a refinement of
Corollary~\ref{AutVbar}. Thus throughout this subsection we assume that $\rho$ is
properly quasi-minuscule.  We have identified
$\hat T_0$   with $\widetilde {\mathcal M}/W_0$. Let $D = \overline T_1\cap T_0$ and 
let $\hat D$ be the preimage of $D$ in
$\hat T_0$. The preimage of
$\hat D$ in
$\widetilde {\mathcal M}$ is the kernel of the weight $\varpi$, and $\hat{D} =
f_\varpi^{-1}(p_0)$ as reduced divisors. If $\Lambda_0\subseteq \Lambda$ is the kernel
of
$\varpi$, then $W_0$ acts on $E_{\rm reg} \otimes \Lambda_0$ and $\hat D = (E_{\rm reg}
\otimes \Lambda_0)/W_0 \subseteq \hat T_0$.

\begin{lemma}\label{hatDbirat} Let $X=\pi_1(D)\cap \pi_1((T_0)_{\rm sing})$. Then $X$
has codimension at least two in $\mathcal{M}$. Moreover, the morphism $\hat D \to D$ is
birational and identifies $\hat D$ with the normalization of $D$.
\end{lemma}
\begin{proof}  A weight $\mu$ induces a homomorphism
$\widetilde{\mathcal{M}} \to E_{\rm reg}$, which we shall just write as $\mu$ again.  By
\cite[Lemma 2.3.5]{Spectral}, there is an open dense subset of the preimage of
$(T_0)_{\rm sing}$ in $\widetilde{\mathcal{M}}$ consisting of points  $x$ such that
$\mu_1(x) =
\mu_2(x)$ for two distinct weights $\mu_1, \mu_2$ of $\rho$ and such that no root
$\alpha$ is the identity on $x$. The inverse image of
$D$ in
$\widetilde{\mathcal{M}}$ is the set of $x$ such that some nonzero weight of $\rho$ is
the identity  on
$x$.  Since the nonzero weights of $\rho$ are roots, 
the inverse image of $\pi_1(D)\cap \pi_1((T_0)_{\rm sing})$ in $\widetilde{\mathcal{M}}$
has codimension at least two in
$\widetilde{\mathcal{M}}$, and since $\widetilde{\mathcal{M}}\to \mathcal{M}$ is
finite, the same is true for $X$. The final statement is then clear.
\end{proof}

Let $\hat T$ be the scheme   $\hat{T}_0 \amalg\overline T_1/\sim$, where $\sim$ is the
identification of a point in $\hat D$ with its image in $D$. Thus the normalization of
$\hat T$ is  $\hat{T}_0 \amalg\overline T_1$, and the morphism $\hat{T}_0 \to T_0$
induces a morphism
$\hat T \to \overline{T}$, and hence an algebra homomorphism $\scrO_{\overline{T}} \to
\scrO_{\hat{T}}$.

\begin{proposition}\label{hatT}
The sheaf $\mathcal{Q}$ is a sheaf of $\scrO_{\hat T}$-modules, compatibly with its
structure as a sheaf of 
$\scrO_{\overline T}$-modules.
\end{proposition}
\begin{proof}  The sheaf
$\pi_1{}_*Hom (\mathcal{V},
\mathcal{V})$ is torsion free and has the Hartogs property.  Thus, if there is an
algebra homomorphism
$\hat\nu_*\scrO_{\hat T}\to\pi_1{}_*Hom (\mathcal{V},
\mathcal{V})$ on the complement $U$ of a closed subset of codimension at
least two, then it extends to a homomorphism over $\hat T$, by taking the composition
$$\hat\nu_*\scrO_{\hat T} \to j_*j^*\hat\nu_*\scrO_{\hat T} \to j_*j^*\pi_1{}_*Hom
(\mathcal{V}, \mathcal{V}) \cong \pi_1{}_*Hom (\mathcal{V}, \mathcal{V}),$$
where $j\colon U\to \mathcal{M}$ is the inclusion. 

Let $X=\pi_1(D)\cap \pi_1((T_0)_{\rm sing})$ and let $U=\mathcal{M}-X$.  By
Lemma~\ref{hatDbirat}, $X$ is a  codimension two subset of $\mathcal{M}$. Set
$U_1=
\mathcal{M}-\pi_1((T_0)_{\rm sing})$ and $U_2 = \mathcal{M}-\pi_1(D)$. Then $U = U_1\cup
U_2$, and it will suffice to check that  $\mathcal{Q}$ is a sheaf of $\scrO_{\hat
T}$-modules over both
$U_1$ and $U_2$. Over $U_1$, $\hat T = \overline{T}$ and we conclude by
Proposition~\ref{ovT}. Over $U_2$, $\overline{T}$ is a disjoint union of
$\overline{T}_1$ and
$T_0$. Hence $\mathcal{Q}|U_2$ is a direct sum $\mathcal{Q}_1|U_2\oplus 
\overline{\mathcal{Q}}|U_2$, and we conclude by Proposition~\ref{T1}.
\end{proof}

We define the sheaf $\mathcal{B}_{\hat{D}}$ to be the subsheaf of
$\mathbb{G}_{m,\hat{T}_0}$ such that $$\mathcal{B}_{\hat{D}}(U) = \{f\in \Gamma(U,
\scrO_U^*): f|U\times _{\hat{T}_0}\hat D =1\}.$$ Thus $\mathcal{B}_{\hat{D}}$ is the 
kernel of the restriction homomorphism
$\mathbb{G}_{m,\hat{T}_0}\to i_*\mathbb{G}_{m,\hat{D}}$, where $i$ is the inclusion of
$\hat D$ in $\hat{T}_0$ and there is a defining exact sequence
$$\{1\} \to \mathcal{B}_{\hat{D}} \to \mathbb{G}_{m,\hat{T}_0}\to
i_*\mathbb{G}_{m,\hat{D}} \to \{1\}.$$

Clearly, $\hat \nu_*\mathcal{B}_{\hat D}$ has the Hartogs property.

\begin{lemma}
The sheaf $\hat \nu_*\mathcal{B}_{\hat D}$ is naturally a subsheaf of
$\pi_1{}_*\Aut({\mathcal V})$.
\end{lemma}

\begin{proof}
Since $T_0$ and $\overline T_1$ meet generically transversely in $D$,
an invertible function on $\hat T_0$ which is $1$ along $\hat D$
 extends to a function on $\hat T$ by defining
it to be $1$ on $\overline T_1$. Thus inclusion induces
an isomorphism from $\hat \nu_*\mathcal{B}_{\hat D}$ to the
sheaf of functions on $\hat T$ that are identically $1$ on
$\overline T_1$, and hence there is an inclusion $\hat \nu_*\mathcal{B}_{\hat D} \to
\nu_*\mathbb{G}_{m,\hat T}$. Since by Proposition~\ref{hatT} the sheaf $\mathcal{Q}$  
is a sheaf of
$\scrO_{\hat T}$-modules, there is a homomorphism $\nu_*\mathbb{G}_{m,\hat T} \to
\pi_1{}_*\Aut({\mathcal V})$ and hence a homomorphism $\hat \nu_*\mathcal{B}_{\hat D}
\to \pi_1{}_*\Aut({\mathcal V})$. It is easy to check that this homomorphism is
injective.
\end{proof}

\begin{lemma}\label{Autprime} Let $\Aut'\mathcal{V}$ be the sheaf of automorphisms of
$\mathcal{V}$ which are the identity on the subbundle $\pi_1{}^*\pi_1{}_*\mathcal{V}$,
so that there is a natural homomorphism $\Aut'\mathcal{V} \to
\Aut\overline{\mathcal{V}}$. Then there is an isomorphism 
$\hat \nu_*\mathcal{B}_{\hat D}\to \pi_1{}_*\Aut'\mathcal{V}$ such that the composition
$$\hat \nu_*\mathcal{B}_{\hat D}\to \pi_1{}_*\Aut'\mathcal{V} \to
\pi_1{}_*\Aut\overline{\mathcal{V}} = \hat\nu_*\mathbb{G}_{m, \hat T_0}$$
is the inclusion.
\end{lemma}
\begin{proof} The previous lemma established an inclusion of $\hat
\nu_*\mathcal{B}_{\hat D}$ in $\pi_1{}_*\Aut({\mathcal V})$. Since the image of $\hat
\nu_*\mathcal{B}_{\hat D}$ is contained in the sheaf of functions on $\hat T$ that are
identically
$1$ on $\overline T_1$,  the  image  of $\hat
\nu_*\mathcal{B}_{\hat D}$ in $\pi_1{}_*\Aut({\mathcal V})$ is
contained in $\pi_1{}_*\Aut'\mathcal{V}$. To find a homomorphism in the other
direction, we use:

\begin{claim} The image of $\pi_1{}_*\Aut'\mathcal{V}$ in 
$\pi_1{}_*\Aut\overline{\mathcal{V}} \cong \hat\nu_*\mathbb{G}_{m,\hat{T}_0}$ is
contained in $\hat \nu_*\mathcal{B}_{\hat D}$.
\end{claim}
\begin{proof} Since $\hat \nu_*\mathcal{B}_{\hat D}$ has the Hartogs property and
$\hat\nu_*\mathbb{G}_{m,\hat{T}_0}$ is torsion free, it suffices to establish the
containment on the complement of a codimension two subset of
$\mathcal{M}$. As before, let
$X=\pi_1(D)\cap \pi_1((T_0)_{\rm sing})$, let $ U  =\mathcal{M}-   X$ let $U_1=
\mathcal{M}-\pi_1((T_0)_{\rm sing})$ and $U_2 = \mathcal{M}-\pi_1(D)$ so that
$U=U_1\cup U_2$. By  Lemma~\ref{hatDbirat}, $X$ is a  codimension two subset of
$\mathcal{M}$ and it will suffice to prove the claim in   $U$. On $U_2$, the
claim is obvious since $\hat\nu_*\mathbb{G}_{m,\hat{T}_0}|U_2=\hat
\nu_*\mathcal{B}_{\hat D}|U_2$. So we may assume that $x\in U_1$. Let $\Omega$ be an
\'etale open neighborhood of $x$ in $U_1$, let $s\in \pi_1{}_*\Aut'\mathcal{V}(\Omega)$,
and let
$f_s$ be the corresponding function on $\Omega\times_{\mathcal{M}}\hat T_0$. If
$\tilde x\in
\Omega\times_{\mathcal{M}}\hat D$ lies over $x$, we must show that $f_s(\tilde x) =1$.

For
$x\in U_1$, let $V_x$ be the restriction of $\mathcal{V}$ to the slice $\{x\}\times
E$, and similarly for $\overline{V}_x$. By the general theory of \cite{FMW},
$\overline{V}_x$ is a regular bundle and hence, if $\nu^{-1}(x) =\{y_1, \dots, y_n\}$
where the $y_i$ are distinct, then  
$\overline{V}_x \cong \bigoplus _{i=1}^n(I_{n_i}\otimes \lambda_i)=\bigoplus
_{i=1}^n\overline{V}_{y_i}$, where the
$\lambda_i$ are line bundles on $E$, pairwise distinct, and indexed by the $y_i$, and
$I_n$ is the unique indecomposable rank
$n$ vector bundle on $E$ with a filtration whose successive quotients are all
isomorphic to $\scrO_E$. Correspondingly, the semilocal ring $H^0(\nu^{-1}(x);
\scrO_{T_0}|\nu^{-1}(x))$ is a direct sum of local rings $\bigoplus_iR_i$, where 
$R_i\cong \Cee[t]/(t^{n_i})$ is the local ring of the fiber $\nu^{-1}(x)$ at $y_i$. Thus
the group of units
$R_i^*$ of
$R_i$ is canonically isomorphic to
$\Cee^*\times S_i$, where $S_i$ is unipotent and the projection $R_i^*\to \Cee^*$ is the
evaluation of an invertible  function on $\Spec R_i$ at the closed point. The ring
$R_i$ is canonically identified with
$\Hom(\overline{V}_{y_i}, \overline{V}_{y_i})$. Using this identification, the
filtration on
$\overline{V}_{y_i}$ is defined by $t^k\overline{V}_{y_i}$ and the action of $z\in
\Cee^*$ is given by multiplication by $z$. 

By definition, $y_i\in \hat D$ if and only if $\lambda_i=\scrO_E$. In this case, by the
local calculations of
\cite[Theorem 2.4.3]{Spectral}, at a generic point of
$D$, the summand of $V_x$ corresponding to the trivial line bundle $\scrO_E$ is of the
form
$\scrO_E^{m_0-1}\oplus I_3$, and there is a natural exact sequence
$$0\to \scrO_E^{m_0}\to \scrO_E^{m_0-1}\oplus I_3 \to I_2 \to 0,$$
such that the quotient $I_2$ is the summand $\overline{V}_{y_i}$ of $\overline{V}_x$.
Let $a$ be the automorphism of $V_x$ induced by the section $s\in
\pi_1{}_*\Aut'\mathcal{V}(\Omega)$ and the choice of the element $\tilde x\in \Omega$
lying over
$x$. By definition, 
$a$ preserves the subbundle
$\scrO_E^{m_0}$ and is the identity on it, and we need to show that $a$ acts trivially
on the associated graded of the quotient
$I_2$. But we can write
$a = \Id + b_1 + b_2$, where $b_1\colon I_3\to I_3$ is multiplication by a scalar $z$
and  $b_3$ is nilpotent, and the induced action on the associated graded of $I_2$ is
multiplication by $1+z$. Restricting to $\scrO_E^{m_0}$, we see that $z=0$, and hence
$a$ induces the constant function $1$, as claimed. Thus, $f_s$ is identically $1$ on
$\Omega\times_{\mathcal{M}}\hat D$.
\end{proof}

Returning to the proof of Lemma~\ref{Autprime}, we have constructed homomorphisms in
both directions between $\pi_1{}_*\Aut'\mathcal{V}$ and $\hat \nu_*\mathcal{B}_{\hat
D}$ which are inverse to each other on a dense open set. Thus, since both sheaves are
torsion free, these homomorphisms induce isomorphisms between 
$\pi_1{}_*\Aut'\mathcal{V}$ and $\hat \nu_*\mathcal{B}_{\hat D}$.
\end{proof}

\begin{corollary}\label{AbarinB} The composition $\overline{\mathcal{A}}_\rho \to
\pi_1{}_*\Aut'\mathcal{V} \to \pi_1{}_*\Aut\overline{\mathcal{V}} \cong \hat \nu_*
\mathbb{G}_{m,\hat T_0}$ defines an embedding of
$\overline{\mathcal{A}}_\rho$ as a closed subsheaf of
$\hat
\nu_*\mathcal{B}_{\hat D}$.
\end{corollary}
\begin{proof} This follows from Proposition~\ref{closure} and the previous result,
using the fact that $\overline{\mathcal{A}}_\rho$ is closed in
$\pi_1{}_*\Aut\mathcal{V}$ and hence in $\pi_1{}_*\Aut'\mathcal{V}$.
\end{proof}

\subsection{Some auxiliary cohomology computations}

We collect some of the general cohomological results needed to compute
$H^0(\mathcal{M}; \mathcal{A})$ and $H^1(\mathcal{M}; \mathcal{A})$. These computations
will ultimately reduce to calculations for
$H^i(\mathcal{M};\hat\nu_*\mathbb{G}_{m,\hat{T}_0})$ and $H^i(\mathcal{M};\hat
\nu_*\mathcal{B}_{\hat D})$, $i=0,1$, or related sheaves, and so we analyze these
cohomology groups here.

\begin{lemma}\label{H^*noD}
Suppose that $\rho$ is minuscule. Then 
$$H^0(\mathcal{M};\hat\nu_*\mathbb{G}_{m,\hat{T}_0}) =\begin{cases} \Cee^*, &\text{if
$E$ is cuspidal;}\\
\Cee^*\times \Zee, &\text{if $E$ is nodal.}
\end{cases}$$
In both cases, $H^1(\mathcal{M};\hat\nu_*\mathbb{G}_{m,\hat{T}_0}) = \{1\}$.
\end{lemma}

\begin{proof} Since $\hat \nu$ is finite, $R^i\hat\nu_*\mathcal{S} = 0$ for every sheaf
$\mathcal{S}$ on $\hat{T}_0$, by e.g.\  \cite[Prop.\ (3.6), p.\ 24]{SGA4half} or
\cite[Corollary 3.4, p.\ 32]{FK}. Thus, by the Leray spectral sequence,
$H^i(\mathcal{M};\hat\nu_*\mathbb{G}_{m,\hat{T}_0})= H^i(\hat{T}_0;
\mathbb{G}_{m,\hat{T}_0})$ for all $i$. For $i=0$, $H^0(\hat{T}_0;
\mathbb{G}_{m,\hat{T}_0})$ is the group of morphisms from $\hat{T}_0$ to $\Gm$, and for
$i=1$, $H^1(\hat{T}_0; \mathbb{G}_{m,\hat{T}_0})=\Pic \hat{T}_0$.

If
$E$ is cuspidal, 
$\hat T_0= \frak h/W_0$. Since $W_0$ is generated by the reflections in the coroots
annihilated by $\varpi$, $\hat T_0\cong \Aff^r$. The result is
immediate in this case.

Now suppose that $E$ is nodal.  Since $\rho$ is minuscule, $\varpi$
is a fundamental weight  and hence   annihilates all
the simple roots except for one, say $\alpha$. Thus, $W_0$ is the Weyl group of a
subroot system $R_0$ of $R$ which is
generated by all the simple roots of $R$ except $\alpha$.  Clearly $\Lambda_0$, the
kernel of $\varpi$ on
$\Lambda$, is exactly the coroot lattice for $R_0$. It follows that $\mathcal{M}_0 =
(E_{\rm reg}\otimes \Lambda_0)/W_0\cong \Aff^{r-1}$. By Lemma~\ref{hatT1}, after an
\'etale base change, the morphism $\hat T_0 \to E_{\rm reg}$ becomes a product
$\mathcal{M}_0 \times C$, where $C\cong E_{\rm reg}$ is a finite \'etale cover of
$E_{\rm reg}$. By Lemma~\ref{geomlemma},  $H^0(\hat{T}_0;
\mathbb{G}_{m,\hat{T}_0})=\Cee^*\times \mathsf{X}(\Gm) = \Cee^*\times \Zee$ and 
$H^1(\hat{T}_0;
\mathbb{G}_{m,\hat{T}_0})=\Pic \hat{T}_0=\{1\}$.
\end{proof}

\begin{lemma}\label{H^*withD}
Suppose that  $\rho$ is properly quasi-minuscule. Then
$$H^0(\mathcal{M}; \hat\nu_*\mathcal{B}_{\hat D})=\begin{cases} \{1\}, &\text{if $E$
is cuspidal;}\\
\Zee =\langle f_\varpi\rangle, &\text{if $E$ is nodal.}
\end{cases}$$ 
In both cases $H^1(\mathcal{M};
\hat\nu_*\mathcal{B}_{\hat D})=\{1\}$.
\end{lemma}

\begin{proof} As in the proof of Lemma~\ref{H^*noD}, an argument using the Leray
spectral se\-quence shows that it suffices to compute
$H^i(\hat{T}_0;\mathcal{B}_{\hat D})$.  We have a short exact sequence of sheaves
$$\{1\} \to \mathcal{B}_{\hat{D}} \to \mathbb{G}_{m,\hat{T}_0}\to
i_*\mathbb{G}_{m,\hat{D}} \to \{1\}.$$

Again we consider
the case when $E$ is cuspidal first. In this case $\hat T_0$ is
again of the form ${\frak h}/W_0$ where $W_0$  is a group generated by
reflections.  Thus $\hat T_0$ is isomorphic to $\Aff^r$. Moreover, $\hat
D=f_{\varpi}^{-1}(1) \cong \mathcal{M}_0$ is  isomorphic to $\Aff^{r-1}$. Thus
$H^0(\hat D;
\Gm) \cong \Cee^*$. Applying the long exact cohomology sequence to the defining exact
sequence for $\mathcal{B}_{\hat{D}}$ then establishes the lemma in this case.

Now suppose that $E$ is nodal. Write $\hat{T}_0 = H/W_0$ and $\hat D =H_0/W_0$, where
$H_0 = \Gm \otimes \Lambda_0 \subseteq H = \Gm\otimes \Lambda$, and the embedding of
$\hat D$ in $\hat{T}_0$ is induced by the inclusion $H_0\subseteq H$. Clearly,
$$H^0(H/W_0; \Gm) =H^0(H;\Gm)^{W_0}= (\Cee^*\times \mathsf{X}(H))^{W_0} =
\Cee^*\times \mathsf{X}(H)^{W_0}.$$
Similarly,  $H^0(H_0/W_0; \Gm)\cong \Cee^* \times \mathsf{X}(H_0)^{W_0}$, and the
restriction map is the identity on the $\Cee^*$ factors and is induced by restriction
of characters on the first factor. We claim that the homomorphism
$\mathsf{X}(H)^{W_0} \to \mathsf{X}(H_0)^{W_0}$ is surjective with kernel $\langle
f_\varpi\rangle$. The exact sequence of tori
$$\{1\} \to H_0 \to H \xrightarrow{\varpi}  \Gm \to \{1\}$$
induces an exact sequence on character groups
$$\{1\} \to \mathsf{X}(\Gm) = \Zee\to \mathsf{X}(H) \to \mathsf{X}(H_0) \to \{1\},$$
where the first inclusion is via $\langle \varpi\rangle$. Taking invariants under
$W_0$, we obtain the sequence
$$\{1\} \to \langle f_\varpi\rangle \to \mathsf{X}(H)^{W_0} \to \mathsf{X}(H_0)^{W_0}
\to H^1(W_0; \Zee),$$
where $H^1(W_0; \Zee)$ is the group cohomology of $W_0$ acting trivially on $\Zee$.
Since every homomorphism from $W_0$ to $\Zee$ is zero, $\mathsf{X}(H)^{W_0} \to 
\mathsf{X}(H_0)^{W_0}$ is surjective with kernel $\langle f_\varpi\rangle$ as claimed.
The result now follows by considering the long exact cohomology sequence associated to
the defining exact sequence for $\mathcal{B}_{\hat{D}}$ as before.
\end{proof}

We will need to consider one more computation related to the existence of an
involution on $\hat T_0$. Suppose that the weights of $\rho$ are invariant under
multiplication by $-1$. Equivalently, $\rho$ is isomorphic to its dual, and it follows
that there is a $G$-invariant nondegenerate bilinear form on $V$, necessarily either
symmetric or skew-symmetric.  In this case, multiplication by
$-1$ on
$E_{\rm reg}$ induces an involution of $\mathcal{M}\times E_{\rm reg}$ which leaves
$T_0$ invariant and lifts to the normalization $\hat{T}_0$. We denote the corresponding
involution of $\hat{T}_0$ by $\tau$. In terms of the Weyl group, let $W_1$ be the
stabilizer of the pair $\{\pm\varpi\}$ and let $w\in W_1$ be an element such that
$w(\varpi) =-\varpi$. Then $w$ acts on $\hat{T}_0$ as the involution $\tau$. Let $\hat
S =\hat{T}_0/\langle \tau\rangle= \widetilde{\mathcal{M}}/W_1$ and let $\varphi\colon
\hat S \to \mathcal{M}$ be the induced morphism. Finally, we set
$$(\hat\nu_*\mathbb{G}_{m,\hat{T}_0})^- = \{f\in \hat\nu_*\mathbb{G}_{m,\hat{T}_0}:
\tau^*f=f^{-1}\}$$
in the appropriate sense. Equivalently, we have the norm homomorphism $N\colon
\hat\nu_*\mathbb{G}_{m,\hat{T}_0} \to \varphi_*\mathbb{G}_{m,\hat{S}}$, and by
definition $(\hat\nu_*\mathbb{G}_{m,\hat{T}_0})^- = \Ker N$. If $f$ is the pullback of
a function $g$ on $\hat S$, $N(f) = g^2$  so that the image of $N$ contains the
image of $\varphi_*\mathbb{G}_{m,\hat{S}}$ under the squaring homomorphism. Thus
$N$ is surjective in the \'etale topology, and there is an exact sequence
$$\{1\} \to (\hat\nu_*\mathbb{G}_{m,\hat{T}_0})^- \to \hat\nu_*\mathbb{G}_{m,\hat{T}_0}
\to \varphi_*\mathbb{G}_{m,\hat{S}} \to \{1\}.$$
Clearly, $(\hat\nu_*\mathbb{G}_{m,\hat{T}_0})^- $ is
a closed subsheaf of $\hat\nu_*\mathbb{G}_{m,\hat{T}_0}$. 

If $\rho$ is properly
quasi-minuscule, then the weights of $\rho$ are always invariant
under $-1$.  Clearly
$w(f_\varpi) = f_\varpi^{-1}$, and so $\hat D = f_\varpi^{-1}(1)$ is invariant under
$\tau$. In fact, since $\varpi$ is a root, the reflection $r_\varpi$ in $\varpi$ is an
element of
$W$ sending $\varpi$ to $-\varpi$ and fixing $\hat D = ( E_{\rm reg}\otimes(\Ker
\varpi))/W_0$ pointwise. Let
$(\hat\nu_*\mathcal{B}_{\hat D})^-= \hat\nu_*\mathcal{B}_{\hat
D} \cap (\hat\nu_*\mathbb{G}_{m,\hat{T}_0})^-$. Since
$\tau|\hat D$ is trivial, a function $f$ on $\hat D$ such that $\tau^*f = f^{-1}$ is
necessarily $\pm 1$. Since there is also an inclusion of $\{\pm 1\}$ in
$(\hat\nu_*\mathbb{G}_{m,\hat{T}_0})^-$, it follows that 
$$(\hat\nu_*\mathbb{G}_{m,\hat{T}_0})^- \cong  (\hat\nu_*\mathcal{B}_{\hat D})^- 
\times \{\pm 1\}.$$ 
Since $(\hat\nu_*\mathbb{G}_{m,\hat{T}_0})^- $ is
a closed subsheaf of $\hat\nu_*\mathbb{G}_{m,\hat{T}_0}$, $(\hat\nu_*\mathcal{B}_{\hat
D})^-$ is a closed subsheaf of $\hat\nu_*\mathcal{B}_{\hat D}$.

\begin{lemma}\label{tau}
In the above notation,
$$H^0(\mathcal{M}; (\hat\nu_*\mathbb{G}_{m,\hat{T}_0})^-) =\begin{cases} \{\pm1\},
&\text{if
$E$ is cuspidal;}\\ \{\pm1\}\times \Zee, &\text{if $E$ is nodal.}\end{cases}$$
Moreover, in both cases $H^1(\mathcal{M}; (\hat\nu_*\mathbb{G}_{m,\hat{T}_0})^-)
=\{1\}$.
\end{lemma}

\begin{proof}
Using the defining exact sequence for $(\hat\nu_*\mathbb{G}_{m,\hat{T}_0})^-$ and 
Lemma~\ref{H^*noD}, it suffices to compute the kernel and cokernel of the
norm homomorphism 
$$H^0(\mathcal{M};\hat\nu_*\mathbb{G}_{m,\hat{T}_0})\to
H^0(\mathcal{M};\varphi_*\mathbb{G}_{m,\hat{S}}).$$
Clearly $H^0(\mathcal{M};\varphi_*\mathbb{G}_{m,\hat{S}}) =
H^0(\hat S;\mathbb{G}_{m,\hat{S}})$ is the group of $\tau$-invariant morphisms  from
$\hat{T}_0$ to $\Gm$. Thus, if $E$ is cuspidal, the norm map on $H^0$ is just
$\Cee^*\to \Cee^*$ given by squaring, so that its kernel is $\{\pm1\}$ and its
cokernel is trivial. In the nodal case, $H^0(\hat{T}_0; \Gm) \cong \Cee^*\times
\Zee$, where we can take $f_\varpi$ to be a generator for $\Zee$. The involution $\tau$
acts trivially on the
$\Cee^*$ factor and via
$-1$ on the
$\Zee$ factor, so that $H^0(\hat S;\mathbb{G}_{m,\hat{S}}) \cong \Cee^*$ and the norm
map $\Cee^*\times \Zee$ is squaring on the first factor and zero on the second. Thus,
in this case, the kernel is $\{\pm1\} \times \Zee$ and the cokernel is trivial.
\end{proof}

Using the isomorphism $(\hat\nu_*\mathbb{G}_{m,\hat{T}_0})^- \cong  (\hat\nu_*\mathcal{B}_{\hat D})^- 
\times \{\pm 1\}$ gives:

\begin{lemma}\label{tau+1}
Suppose that $\rho$ is properly quasi-minuscule. Then the
weights of $\rho$ are invariant under $-1$ and $\hat
D\subset \hat T_0$ is non-empty and is invariant under $-1$. Moreover
$$H^0(\mathcal{M}; (\hat\nu_*\mathcal{B}_{\hat D})^-) =\begin{cases} \{1\},
&\text{if
$E$ is cuspidal;}\\ \Zee =\langle f_\varpi\rangle, &\text{if $E$ is nodal.}\end{cases}$$
In both cases $H^1(\mathcal{M}; (\hat\nu_*\mathcal{B}_{\hat D})^-)
=\{1\}$. \qed
\end{lemma}

\section{Explicit examples}

In this section, we give a concrete description of the sheaf
$\overline{\mathcal{A}}_\rho$ in case $\rho$ is the standard representation of $SL(n+1,
\Cee)$, $Sp(2n, \Cee)$, or
$Spin(n, \Cee)$ or one of the minuscule representations of a group of type $E_6$ or
$E_7$. These cases cover every group which has a minuscule representation, and all of
the minuscule representations with the exception of the higher exterior powers of the
standard representation of $SL(n+1, \Cee)$ and the spin or half spin representations
for $Spin(n, \Cee)$. In the next section, we will give a general procedure for dealing
with all of the properly quasi-minuscule  representations and thus give a
uniform proof of Thereom~\ref{main} for all groups.

We consider the various cases individually.
\subsection{Type $A_n$: $G=SL(n+1, \Cee)$}

Let $\rho\colon SL(n+1, \Cee) \to GL(n+1, \Cee)$ be the defining
(minuscule) representation and let $\mathcal{V}$ be the corresponding vector bundle
over $\mathcal{M}\times E$. The   restriction
$\mathcal V|(\{x\}\times E)$ is regular for every $x\in{\mathcal M}$. It follows that
$T=T_0= \hat T_0$. Since $\rho$ is faithful, $\overline {\mathcal A}_\rho ={\mathcal
A}$.

\begin{proposition}
${\mathcal A}\subseteq \nu_*\mathbb{G}_{m,T}$ is the kernel
of the norm map
$N\colon \nu_*\mathbb{G}_{m,T}\to \mathbb{G}_{m,{\mathcal M}}$.
\end{proposition}

\begin{proof}
This is immediate from the identification of $N$ with the determinant homomorphism.
\end{proof}

\begin{corollary}
There is an exact sequence
$$\{1\}\to {\mathcal A}\to \nu_*\mathbb{G}_{m,T}\xrightarrow{N} \mathbb{G}_{m,{\mathcal
M}}\to
\{1\}.$$
\end{corollary}
\begin{proof} It suffices to observe that the image of $N$ contains the
$(n+1)^{\textrm{st}}$ powers of functions in $\mathbb{G}_{m,{\mathcal
M}}$ and hence $N$ is surjective in the \'etale topology.
\end{proof}

Thus, there is a  exact sequence 
$$0\to H^0(\mathcal{M}; {\mathcal A})\to H^0(T; \mathbb{G}_{m,T}) \to
H^0(\mathcal{M};
\mathbb{G}_{m,{\mathcal{M}}}) \to  H^1(\mathcal{M}; {\mathcal A})\to \Pic T \to \Pic
\mathcal{M} \to 0.$$ 
Since
${\mathcal M}\cong \Aff^n$, 
$H^0({\mathcal M};\mathbb{G}_{m, {\mathcal M}})=\Cee^*$. By Lemma~\ref{H^*noD} (or
directly),
$\Pic T =0$, and $H^0(T; \mathbb{G}_{m,T}) \cong\Cee^*$, in the cuspidal case,  and
$H^0(T;
\mathbb{G}_{m,T}) \cong\Cee^*\times \Zee$ in the nodal case. The norm map $N$ is raising
to the $(n+1)^{\textrm{st}}$ power on the
$\Cee^*$ factors and is trivial on the $\Zee$ factor, if it is present (because
$N(f_\varpi) =\det$). Hence
$N$ is surjective, so that $H^1(\mathcal{M}; {\mathcal A}) =\{1\}$, and $
H^0(\mathcal{M}; {\mathcal A}) = \Ker N$ is either
$\Zee/(n+1)\Zee$ or $(\Zee/(n+1)\Zee)\times
\Zee$, depending on whether $E$ is cuspidal or nodal.

\begin{remark} More generally, let $E$ be any Weierstrass cubic, let $B$ be an arbitrary
scheme, and let
$\mathcal{V}\to B\times E$ be a family of regular semistable bundles of rank $n+1$ with
trivial determinant. In case $E$ is singular, we do not require that 
$\mathcal{V}|\{x\}\times E$ pulls back to the trivial vector bundle on $\widetilde E$
for every $x\in B$.  Let
$\nu\colon T_{\mathcal{V}} \to B$ be the spectral cover; it is compatible with base
change. In this case, if
$\Xi_{\mathcal{V}}$ is the associated $SL(n+1, \Cee)$-bundle, then arguments as above
show that there is an exact sequence
$$\{1\} \to \pi_1{}_*\Aut
\Xi_{\mathcal{V}}\to\nu_*\mathbb{G}_{m,T_{\mathcal{V}}}\xrightarrow{N}
\mathbb{G}_{m,B}\to
\{1\}.$$
A similar remark holds for elliptic fibrations. 
\end{remark}

\subsection{Type $C_n$: $G=Sp(2n,\Cee)$}

Let $\rho\colon G \to GL(2n, \Cee)$ be the standard representation. Thus $\rho$ is
injective and minuscule. In particular $\overline{\mathcal{A}}_\rho =\mathcal{A}$. The
corresponding vector bundle
$\mathcal{V}$ has regular restriction to each slice $\{x\}\times E$ since the
difference of any two distinct weights is a root. Hence
$T=T_0=\hat T_0$ in this case.  

There is a nondegenerate symplectic form $B$ on $\mathcal{V}$,  and
hence an isomorphism
$\varphi \colon \mathcal{V}\to \mathcal{V}^*$ via $\varphi(v)(w) = B(v,w)$. Note that
$A$ is symplectic if and only if $\varphi\circ A^{-1} = A^*\circ \varphi$. As we have
seen, the isomorphism $\varphi$ induces an involution $\tau$ on
$T=\hat T_0$ covering the identity on $\mathcal{M}$. We can see $\tau$ quite
concretely as follows: we have $T = \mathbf{Spec} \,\pi_1{}_*Hom(\mathcal{V},
\mathcal{V})$, viewed as a sheaf of $\scrO_{\mathcal{M}}$-algebras, and the involution
is the one defined by
$$\tau (A) = \varphi^{-1}\circ A^* \circ \varphi.$$
It is easy to check that this agrees with the involution $\tau$ on $T$ as defined
previously.  In particular, $A$ is symplectic if and only if $\tau(A) = A^{-1}$. By
definition, we have:

\begin{proposition}$\mathcal{A}=(\nu_*\mathbb{G}_{m,T})^-$. \qed
\end{proposition}

Applying Lemma~\ref{tau} then gives:

\begin{corollary} $H^1(\mathcal{M}; {\mathcal
A})=\{1\}$, and $H^0(\mathcal{M}; {\mathcal A})$ is isomorphic to $\{\pm
1\}$ if $E$ is cuspidal and  $\{\pm 1\}\times\Zee$ if $E$ is nodal.
\end{corollary}

As in the case of $SL(n+1, \Cee)$, there is a relative version of the sheaf
$\mathcal{A}$, including the case where $E$ is smooth or certain cases where the
$G$-bundles $\Xi|\{x\}\times E$ need not pull back to the trivial bundle on the
normalization.

\subsection{Type $D_n$: $G=Spin(2n, \Cee)$}

In this case,  let $\rho$ be the $2n$-dimensional orthogonal
representation, i.e.\ $\rho$ is the composition $Spin(2n, \Cee) \to SO(2n, \Cee) \to
GL(2n, \Cee)$. In this case the kernel of
$\rho$ is a subgroup of order
$2$. The representation $\rho$ is minuscule, so that
$T=T_0$. However, since there are pairs of weights whose difference is not a
multiple of a root, there exist
$x\in
\mathcal{M}$ such that  the restriction
${\mathcal V}|(\{x\}\times E)$ is not regular. Thus $T_0$ is not normal,  and hence
$\hat T_0\neq T_0$.

As in the symplectic case, the quadratic form on $\mathcal{V}$ defines an
isomorphism $\varphi \colon \mathcal{V}\to \mathcal{V}^*$. The involution
$\tau\colon \hat T_0 \to
\hat T_0$  is clearly induced by the involution of $\pi_1{}_*Hom (\mathcal{V},
\mathcal{V})$ defined by $A\mapsto \varphi^{-1}\circ A^*\circ \varphi$. Moreover, $A$
is orthogonal if and only if $\tau(A) = A^{-1}$. Thus, via Corollary~\ref{AutVbar},

\begin{proposition} $\overline{\mathcal{A}}_\rho = (\hat \nu_*\mathbb{G}_{m, \hat
T_0})^-$.\qed
\end{proposition}

\begin{corollary} $H^1(\mathcal{M};
{\mathcal A})=\{1\}$, and $H^0(\mathcal{M}; {\mathcal A})$ is isomorphic to $Z(G)$ if
$E$ is cuspidal and $Z(G)\times \Zee$ if $E$ is nodal. 
\end{corollary}
\begin{proof} By Lemma~\ref{tau}, $H^1(\mathcal{M};
\overline{\mathcal{A}}_\rho)=\{1\}$, and $H^0(\mathcal{M};
\overline{\mathcal{A}}_\rho)$ is either $\{\pm1\}$ if $E$ is cuspidal or $\Zee\times
\{\pm1\}$ if $E$ is nodal. By Lemma~\ref{AvsAbar}, $H^1(\mathcal{M}; \mathcal A)  =
\{1\}$ and there is an exact sequence
$$0\to \Zee/2\Zee \to  H^0(\mathcal{M}; \mathcal A) \to H^0(\mathcal{M}; \overline
{\mathcal A}_\rho) \to \{1\}.$$
It follows that the torsion subgroup of $H^0(\mathcal{M}; \mathcal A)$ has order $4$
and its rank is zero if $E$ is cuspidal and one if $E$ is nodal. Since $H^0(\mathcal{M};
\mathcal A)$ contains a subgroup isomorphic to $Z(G)$, it follows that
$H^0(\mathcal{M}; \mathcal A)$ is as claimed.
\end{proof}

\subsection{Type $B_n$: $G=Spin(2n+1, \Cee)$}

In this case, let $\rho$ be the standard representation induced by $Spin(2n+1) \to
SO(2n+1, \Cee) \subseteq GL(2n+1, \Cee)$. The kernel of $\rho$ has order two. 
  The representation $\rho$ is properly quasi-minuscule
and the trivial weight has multiplicity one.  The spectral cover is
$T=T_0\cup T_0$, where
$T_0$ is reduced. It is straightforward to check that, if $E$ is cuspidal, then $T_0$ is
normal, whereas $T_0$ is not normal in case $E$ is nodal. The quadratic form on
$\mathcal{V}$ induces an involution on $T_0$ and $\hat T_0$. With this said, we have:

\begin{proposition} $\overline{\mathcal{A}}_\rho = (\hat\nu_*{\mathcal B}_{\hat D})^-$.
\end{proposition}
\begin{proof} By Corollary~\ref{AbarinB}, $\overline{\mathcal{A}}_\rho$
is a closed subsheaf of $\hat\nu_*{\mathcal B}_{\hat D}$. On the other hand,
$(\hat\nu_*{\mathcal B}_{\hat D})^-$ is also a closed subsheaf of $\hat\nu_*{\mathcal
B}_{\hat D}$. Thus, it suffices to prove that the restrictions of
$\overline{\mathcal{A}}_\rho$ and $(\hat\nu_*{\mathcal B}_{\hat D})^-$ to the \'etale open
set $U$ of Lemma~\ref{Gsplit} and Corollary~\ref{split} agree. On $U$, $\mathcal{V} =
\bigoplus _\mu \mathcal{L}_\mu$, where $\overline{\mathcal{A}}_\rho$ acts trivially on
$\mathcal{L}_0$, acts on
$\mathcal{L}_\mu$ via the character $\mu$, and hence on $\mathcal{L}_{-\mu}$ via the
character
$\mu^{-1}$. Clearly, this identifies $\overline{\mathcal{A}}_\rho|U$ with
$(\hat\nu_*{\mathcal B}_{\hat D})^-|U$.
\end{proof}

Standard arguments, using Lemma~\ref{tau+1}, then show:

\begin{corollary} $H^1(\mathcal{M}; {\mathcal A})=\{1\}$, and 
$H^0(\mathcal{M}; {\mathcal A})$ is isomorphic to $Z(G)$ if $E$ is cuspidal and 
$Z(G)\times \Zee$ if $E$ is nodal.  \qed
\end{corollary}

\subsection{Type $E_6$}

Let $G$ be of type $E_6$ and let $\rho\colon G \to GL(V)$ be one of the $27$-dimensional
minuscule representations of $G$. The representation $\rho$ is injective and 
minuscule. Thus $T=T_0$. However, $T_0$ is never normal, since there are pairs of
weights of $\rho$ whose difference is not a multiple of a root. The group
$G$ is the stabilizer of a symmetric cubic form $C$ in the automorphism group of $V$.
If we decompose $V$ into a direct sum of weight spaces $V=\bigoplus_\mu L_\mu$, then 
the form
$C$ preserves the decomposition of $V$ into weight spaces. It is non-trivial on
$L_{\mu_1}\otimes L_{\mu_2}\otimes L_{\mu_3}$ if and only if $\mu_1+\mu_2+\mu_3=0$. 

\begin{claim}\label{Hprime}  Let $\widetilde{H}$ be the set of  $A\in GL(V)$ such
that
$A|L_\mu$ is multiplication by a nonzero scalar $a_\mu$ and such that the $a_\mu$
satisfy: for every triple of weights $(\mu_1, \mu_2, \mu_3)$ such $\mu_1+\mu_2+\mu_3
=0$,
$a_{\mu_1}a_{\mu_2} a_{\mu_3}=1$. Then $\rho(H)=\widetilde{H}$.
\end{claim}
\begin{proof} Clearly, if $h\in H$, then $\rho(h)|L_\mu$ is multiplication by
$a_\mu=\mu(h)$ and hence the $a_\mu$ satisfy the multiplicative property of the claim.
Conversely, if $A\in \widetilde{H}$, then $A$ preserves $C$, so that
$A=\rho(g)$ for some $g\in G$. As $A$ commutes with $\rho(H)$ and $\rho$ is faithful,
$g$ commutes with
$H$ and hence $g\in H$.
\end{proof}

We turn now to a description of $\mathcal{A}=\overline{\mathcal{A}}_\rho$. By
Corollary~\ref{AutVbar},  $\mathcal{A}$ is a closed subsheaf of
$\hat \nu_*\mathbb{G}_{m, \hat{T}_0}$. We will describe this subsheaf by means of a
correspondence on $\hat T_0$. Recall that $W_0$ is the stabilizer of the highest weight
$\varpi$ of $\rho$. 

\begin{lemma}\label{Weyltrans1} If $\mu_1$ is a weight of $\rho$, then there exist
weights
$\mu_2,
\mu_3$ of $\rho$ such that $\mu_1+\mu_2+\mu_3=0$. The Weyl group
$W$ acts transitively on the set of ordered pairs
$(\{\mu_1, \mu_2,\mu_3\},\mu_1)$, where $\{\mu_1, \mu_2,\mu_3\}$ is an unordered triple
of weights of $\rho$ whose sum is zero.
\end{lemma}
\begin{proof} 
Since the Weyl group acts transitively on the weights of
$\rho$, it suffices to show that the stabilizer $W_0$ of $\varpi$
acts transitively on the set of weights $\mu$ with the
property that $\varpi+\mu$ is the negative of a weight of $\rho$.
The kernel of $\varpi$ in the coroot lattice of $E_6$ is the coroot lattice of $D_5$.
Thus $W_0$ contains a subgroup isomorphic to the Weyl group of $D_5$.  
As a representation of $Spin(10)$,   $V$ is a direct sum $\Cee\oplus
V^{10}\oplus S^{16}$, where $V^{10}$ is the orthogonal representation and $S^{16}$ is
one of the half spinor representations. The weight spaces $L_\mu$, where
$\mu$ has  the property that
$\varpi+\mu$ is the negative of a weight of $\rho$, are all weight spaces in $V^{10}$.
Since the Weyl group of $D_5$ acts transitively on these weights, the result
follows.
\end{proof}

In terms of cubic surfaces, the weights of $\rho$ may be identified with lines on a
smooth cubic surface, and the lemma says that the stabilizer $W_0$ in $W$ of a line
$\ell$ acts transitively on the $10$ lines which meet $\ell$.

Given the weight $\varpi$, fix weights $\mu_2, \mu_3$ such that $\varpi+\mu_2+\mu_3
=0$. Let $W'$ be the stabilizer of the unordered triple $\{\varpi, \mu_2, \mu_3\}$ and
let $W''$ be the stabilizer of the ordered pair $(\{\varpi, \mu_2, \mu_3\}, \varpi)$.
Thus $W''$ is a subgroup both of $W_0$ and $W'$. Set $S'=\widetilde{\mathcal{M}}/W'$
and
$S''=\widetilde{\mathcal{M}}/W''$. Thus there is a diagram of schemes
$$\begin{matrix}
& & S'' && &\\
& \swarrow && \searrow &\\
\hat T_0 & && & S'\\
& \searrow & &\swarrow &\\
& & \mathcal{M} && &
\end{matrix}.$$
Let $\nu'\colon S' \to \mathcal{M}$ and $\nu''\colon S'' \to \mathcal{M}$ be the induced
morphisms. Then there is a homomorphism $c\colon \hat\nu_*\mathbb{G}_{m, \hat T_0} \to
\nu'_*\mathbb{G}_{m, S'}$ defined by pulling a function on $\hat T_0$ back to
$S''$, followed by the norm homomorphism $\nu''_*\mathbb{G}_{m, S''} \to
\nu'_*\mathbb{G}_{m,  S'}$. By Lemma~\ref{Hartogssub}, $\Ker c$ is a closed subsheaf of
$\hat\nu_*\mathbb{G}_{m, \hat T_0}$. If $x\in \mathcal{M}$ is a point such that the
morphism $\widetilde{\mathcal{M}} \to \mathcal{M}$ is \'etale and $\tilde x\in
\widetilde{\mathcal{M}}$ lies in the fiber over $x$,  the fiber
$\nu^{-1}(x)$ is identified with the set of weights of $\rho$ and the fiber 
$(\nu')^{-1}(x)$ is identified with the set of unordered triples $\{\mu_1, \mu_2,
\mu_3\}$ such that 
$\mu_1+\mu_2+\mu_3 =0$. Clearly, 
$$c(f)(\{\mu_1, \mu_2, \mu_3\}) = f(\mu_1)f(\mu_2)f(\mu_3)$$
and hence $c(f) =1$ if and only if $f(\mu_1)f(\mu_2)f(\mu_3)=1$ for all such triples.

\begin{proposition}\label{AforE6} $\mathcal{A} =\Ker c$.
\end{proposition}
\begin{proof} Since both $\mathcal{A}$ and $\Ker c$ and  closed subsheaves of the
torsion free sheaf
$\hat\nu_*\mathbb{G}_{m, \hat T_0}$, it suffices to check the equality on an \'etale
open dense set. Let $U$ be the \'etale open set of Lemma~\ref{Gsplit} and
Corollary~\ref{split}. The result over $U$ follows immediately from the above
description of $\Ker c$ and Claim~\ref{Hprime}.
\end{proof}

\begin{proposition} $H^1(\mathcal{M}; {\mathcal A})=\{1\}$, and $H^0(\mathcal{M};
{\mathcal A})$ is isomorphic to $\Zee/3\Zee$ if $E$ is cuspidal and 
$(\Zee/3\Zee)\times \Zee$ if $E$ is nodal. 
\end{proposition}
\begin{proof} Let $\mathcal{I}= \operatorname{Im}c$. Then $\mathcal{I}$ is a subsheaf
of $\nu'_*\mathbb{G}_{m,  S'}$ and  there is an exact sequence
$$\{1\} \to \mathcal{A} \to \hat\nu_*\mathbb{G}_{m, \hat T_0} \to \mathcal{I}\to
\{1\}.$$
First suppose that $E$ is cuspidal. Since there is a finite morphism from
$\widetilde{\mathcal{M}}=\frak h$ to $S'$, $\Gamma(S', \mathbb{G}_{m,S'}) \subseteq
\Gamma(\frak h, \mathbb{G}_{m,{\frak h}^*})=\Cee^*$ as $\frak h$ is affine. Thus
$\Gamma(S',
\mathbb{G}_{m,S'}) =\Cee^*$ consists of the nonzero constant functions. It follows that
the homomorphism $c\colon H^0(\mathcal{M}; \hat\nu_*\mathbb{G}_{m, \hat T_0} ) \to
H^0(\mathcal{M}; \nu'_*\mathbb{G}_{m,  S'})$ is the cubing homomorphism from $\Cee^*$
to $\Cee^*$ and factors through the subgroup $H^0(\mathcal{M};\mathcal{I})$ of
$H^0(\mathcal{M};
\nu'_*\mathbb{G}_{m,  S'})$. This forces
$H^0(\mathcal{M};\mathcal{I}) = H^0(\mathcal{M}; \nu'_*\mathbb{G}_{m,  S'}) =\Cee^*$,
and the induced homomorphism $$H^0(\mathcal{M}; \hat\nu_*\mathbb{G}_{m, \hat T_0} ) \to
H^0(\mathcal{M};\mathcal{I})$$ is again the cubing homomorphism $\Cee^*\to \Cee^*$. Thus
it is surjective with kernel $\Zee/3\Zee$. It follows that $H^0(\mathcal{M}; {\mathcal
A})=\Zee/3\Zee$ and that $H^1(\mathcal{M}; {\mathcal A}) \subseteq H^1(\mathcal{M};
\hat\nu_*\mathbb{G}_{m, \hat T_0} )$, which is trivial by Lemma~\ref{H^*noD}. This
concludes the proof in the cuspidal case.

Now assume that $E$ is nodal. By Lemma~\ref{H^*noD}, $H^0(\mathcal{M};
\hat\nu_*\mathbb{G}_{m, \hat T_0} ) \cong \Cee^* \times \Zee$. As in the proof of the
cuspidal case, it will suffice to prove that $H^0(\mathcal{M};\mathcal{I}) = 
H^0(\mathcal{M}; \nu'_*\mathbb{G}_{m,  S'}) =\Cee^*$, and that the induced homomorphism
$H^0(\mathcal{M}; \hat\nu_*\mathbb{G}_{m, \hat T_0} ) \to
H^0(\mathcal{M};\mathcal{I})$ is the homomorphism from $\Cee^* \times \Zee$ to $\Cee^*$
given by $(z,n) \mapsto z^3$. Clearly $c(f_\varpi) =1$, so that it suffices to show
that $H^0(\mathcal{M}; \nu'_*\mathbb{G}_{m,  S'})=H^0(S'; \mathbb{G}_{m,S'})=\Cee^*$. As
in the proof of Lemma~\ref{H^*withD}, 
$$H^0(S'; \mathbb{G}_{m,S'}) = H^0(H; \mathbb{G}_{m,H})^{W'} = (\Cee^* \times
\mathsf{X}(H))^{W'} = 
\Cee^* \times \mathsf{X}(H)^{W'}.$$
Since $\mathsf{X}(H) =\Lambda^*$, it suffices to prove that $(\frak
h^*)^{W'} =0$, or equivalently that $\frak h^{W'}=0$. It is easy to check that,
if $\varpi, \mu_2,
\mu_3$ are three weights of $\rho$ whose sum is zero, then $\{\varpi, \mu_2,
\mu_3\}^\perp$ is a root lattice of type $D_4$ and hence $W'$ contains the Weyl group
of $D_4$. Also, by Lemma~\ref{Weyltrans1}, $W'$ acts transitively on the set $\{\varpi,
\mu_2, \mu_3\}$. From this, it follows that there are no $W'$-invariant vectors in
$\frak h =\operatorname{span}\{\varpi, \mu_2, \mu_3\} \oplus \{\varpi, \mu_2,
\mu_3\}^\perp$.

Thus $H^0(S'; \mathbb{G}_{m,S'})=\Cee^*$, and the rest of the proof proceeds as in the
cuspidal case.
\end{proof}

\subsection{Type $E_7$}

Let $G$ be of type $E_7$, and let $\rho\colon G \to GL(V)$ be the $56$-dimensional
irreducible representation of $E_7$.
Since $\rho$ is injective and minuscule, $\overline{\mathcal{A}}_\rho =\mathcal{A}$
and $T=T_0$.  Since there are pairs of weights of $\rho$ whose difference is
not a multiple of a root,
$T_0$ is never normal. The group $G$ is the stabilizer of a non-degenerate
skew form  and a symmetric quartic form on $V$.
Decomposing $V$ as a direct sum of weight spaces, $V=\bigoplus_\mu
L_\mu$, the skew form is non-zero exactly on tensor products of
the form $L_\mu\otimes L_{-\mu}$ and the quartic form is nontrivial exactly on tensor
products of the form $L_{\mu_1}\otimes L_{\mu_2}\otimes L_{\mu_3}\otimes L_{\mu_4}$ when
$\mu_1+\mu_2+\mu_3+\mu_4=0$ and $\mu_1+\mu_i\not=0$ for $i=2,3,4$.
We say that an unordered quadruple $\{\mu_1,\mu_2,\mu_3,\mu_4\}$ of weights
\textsl{properly sums to zero} $\mu_1+\mu_2+\mu_3+\mu_4=0$ and $\mu_i+\mu_j\neq 0$ for
all $i,j$, $1\leq i,j\leq 4$.

An argument as for the case of $E_6$ then shows:

\begin{claim}\label{char}
Let $\widetilde{H}$ be the set of $A\in GL(V)$ such that $A|L_\mu$ is
multiplication by a scalar $a_\mu$ and such that the $a_\mu$
satisfy:
 $a_\mu a_{-\mu}=1$ for all weights $\mu$ and
$a_{\mu_1}a_{\mu_2}a_{\mu_3}a_{\mu_4}=1$ whenever
$\{\mu_1,\mu_2,\mu_3,\mu_4\}$  properly sums to zero.
Then $\rho(H) = \widetilde{H}$.
\qed
\end{claim}

We turn now to the description of $\mathcal{A}$. Since the skew form identifies $V$ and
$V^*$, there is an involution $\tau$ on $\hat T_0$. By Claim~\ref{char}, $\mathcal{A}$
is a closed subsheaf of $(\hat \nu_*\mathbb{G}_{m, \hat T_0})^-$.  As in the case of
$E_6$, we will describe this subsheaf by a correspondence. We begin with a lemma on the
Weyl group:

\begin{lemma}\label{E7trans} The group $W$ acts transitively on the  
set of all  pairs
$(\{\mu_1,\mu_2,\mu_3,\mu_4\},\mu_1)$ consisting of a quadruple which 
properly sums to zero and one of the weights in the quadruple. Moreover, if 
$\{\varpi,\mu_2,\mu_3,\mu_4\}$ is such a quadruple, the
stabilizer
$W_0$ of  $\varpi$ contains an element which cyclically permutes the set
$\{\mu_2,\mu_3,\mu_4\}$.
\end{lemma}

\begin{proof}
Since $W$ acts transitively on the weights of $\rho$,
to prove the first statement, it suffices to show that $W_0$ acts transitively on the
set of all quadruples $\{\varpi,\mu_2,\mu_3,\mu_4\}$ which properly sum to zero. The
kernel of
$\varpi$ is the coroot lattice $\Lambda_1$ of
$E_6$, and hence   $W_0$ contains a subgroup $W_1$ isomorphic to the Weyl group of
$E_6$. Let $G_1$ be the corresponding subgroup of $G$  of type $E_6$. The
representation 
$V$ viewed as a representation of  $G_1$, is isomorphic to a
direct sum
$\Cee\oplus \Cee \oplus A\oplus \overline{A}$ where $A$ and
$\overline A$ are the two non-isomorphic minuscule representations
of $G_1$. The weight spaces $L_{\mu}$ with $\mu$ as above form one
of the minuscule factors $A$ or $\overline A$. In fact, possibly after relabeling, the
quadruples containing 
$\varpi$ which properly sum to zero  consist of $\varpi$ and a triple of weights
which become weights for the maximal torus of $G_1$ acting on $A$ and which sum to zero
as weights acting on $\Lambda_1$. Since $W_1$ acts transitively on these
weights, the first statement follows. The second is a consequence of
Lemma~\ref{Weyltrans1}.
\end{proof}

Fix a quadruple $\{\varpi,\mu_2,\mu_3,\mu_4\}$ which properly sums to zero. Let $W'$ be
the stabilizer of $\{\varpi,\mu_2,\mu_3,\mu_4\}$ and let $W''$ be the stabilizer of the
pair $(\{\varpi,\mu_2,\mu_3,\mu_4\}, \varpi)$. Thus $W''$ is a subgroup both of $W_0$
and of $W'$. As in the case of
$E_6$, we form the spaces $S'=\widetilde{\mathcal{M}}/W'$ and
$S''=\widetilde{\mathcal{M}}/W''$. Clearly, the degree of the morphism $S''\to S'$ is
$4$. If
$\nu'\colon S'\to
\mathcal{M}$ is the induced morphism, then we can again define the homomorphism
$c\colon \hat
\nu_*\mathbb{G}_{m,\hat T_0}\to
\hat \nu'_*\mathbb{G}_{m,S'}$. Clearly, the involution $\tau$ extends to involutions
on $S''$ and $S'$, also denoted $\tau$, and $c$ is equivariant. Thus the restriction of
$c$ to $(\hat \nu_*\mathbb{G}_{m,\hat T_0})^-$ takes values in 
$(\hat \nu'_*\mathbb{G}_{m,S'})^-$. An argument as in the proof of
Proposition~\ref{AforE6} gives:

\begin{proposition} $\mathcal{A} = (\hat \nu_*\mathbb{G}_{m,\hat T_0})^- \cap \Ker c$.
\qed
\end{proposition}

\begin{proposition} $H^1(\mathcal{M};{\mathcal A})=\{1\}$, and $H^0(\mathcal{M};
{\mathcal A})$ is isomorphic to $\Zee/2\Zee$ if $E$ is cuspidal and
$(\Zee/2\Zee)\times \Zee$ if $E$ is nodal. 
\end{proposition}

\begin{proof} Let $\mathcal{I}= \operatorname{Im}(c|(\hat \nu_*\mathbb{G}_{m,\hat
T_0})^-)$. Then $\mathcal{I}$ is a subsheaf of $(\nu'_*\mathbb{G}_{m,
S'})^-$ and  there is an exact sequence
$$\{1\} \to \mathcal{A} \to (\hat\nu_*\mathbb{G}_{m, \hat T_0})^-
 \to \mathcal{I}\to \{1\}.$$
 The proof is completed, using Lemma~\ref{tau},
 by showing that $H^0({\mathcal
 M};{\mathcal I})=0$.

We begin by showing that $H^0({\mathcal M};(\nu'_*\mathbb{G}_{m,S'})^-)=\{\pm1\}$.
 In the case when $E$ is cuspidal, $S'$ is a quotient of ${\frak
 h}$ by a finite group and hence $H^0({\mathcal M};\nu'_*\mathbb{G}_{m,S'})=\Cee^*$.
Since $\tau$ acts trivially on this $\Cee^*$,    $H^0({\mathcal
M};(\nu'_*\mathbb{G}_{m,S'})^-)=\{\pm1\}$. When $E$ is nodal,  $H^0({\mathcal
M};\nu'_*\mathbb{G}_{m,S'})$
 is equal to the product $\Cee^*\times \mathsf{X}(H)^{W'}$ and we must compute the
second factor. Given the quadruple $\{\varpi,\mu_2,\mu_3,\mu_4\}$ which properly sums to
zero, let $\Lambda'$ be the intersection of the kernels of the weights
$\varpi,\mu_2,\mu_3,\mu_4$ and let $H'\subseteq H$ be the corresponding closed
connected subgroup. Note that $\Lambda'$ is a coroot lattice of type $D_4$. Let $H''$ be
the torus whose Lie algebra is spanned by the weights, which we view as elements of
$\frak h$ via the inner product.  Note that $W'$ acts on $H'$ and
$H''$, and there is a
$W'$-equivariant isogeny
$H'\times H'' \to H$. Correspondingly, there is a $W'$-equivariant embedding
$\mathsf{X}(H) \to \mathsf{X}(H') \oplus \mathsf{X}(H'')$. Since $W'$ contains a
subgroup generated by the reflections in the coroots contained in $\Lambda'$, which is
the Weyl group of $D_4$, it follows that $\mathsf{X}(H')^{W'}=0$. As for the second
factor $\mathsf{X}(H'')\cong \Zee^3$, by  Lemma~\ref{E7trans}
there is an element of $W'$ which fixes $\varpi$ and cyclically permutes
$\mu_2,\mu_3,\mu_4$. On the other hand, again by  Lemma~\ref{E7trans}, there is an
element of $W'$ which sends $\varpi$ to $\mu_2$. These two elements generate a subgroup
of the symmetric group $S_4$ which is of index at most two and hence  contains the
alternating group. It follows that $\mathsf{X}(H'')^{W'} =0$ as well. Thus again 
$H^0({\mathcal M};\nu'_*\mathbb{G}_{m,S'})=\Cee^*$ and  
$H^0({\mathcal M};(\nu'_*\mathbb{G}_{m,S'})^-)=\{\pm1\}$.

Thus $H^0({\mathcal M};(\nu'_*\mathbb{G}_{m,S'})^-)=\{\pm1\}$, where the functions in
question are constant functions on $S'$. The final step is to show that the image of
$H^0(\mathcal{M}; \mathcal{I})$ in $H^0({\mathcal
M};(\nu'_*\mathbb{G}_{m,S'})^-)=\{\pm1\}$ is trivial. Let $x_0\in \mathcal{M}$ be the
image of the identity in $\widetilde{\mathcal{M}}$. Since the identity is fixed by $W$,
there is a unique preimage of $x_0$ in $\hat T_0$, as well as in $S'$ and $S''$, and
all of these points are $\tau$-invariant. Let $U\to \mathcal{M}$ be an \'etale
neighborhood of $x_0$ and $u$ a point in the preimage of $x_0$. Then there is a unique
point $\hat u\in U\times_{\mathcal{M}}\hat T_0$ lying above $u$, and similarly there
are unique points $u'\in U\times_{\mathcal{M}}S'$ and $u''\in U\times_{\mathcal{M}}S''$
above $u$. In particular, $\tau(u'') =u''$, where we continue to denote by $\tau$ the
induced involution on $U\times_{\mathcal{M}}S''$. Hence, if $f\in (\hat
\nu_*\mathbb{G}_{m, \hat T_0})^-(U)$, then $f(\hat u)=\pm 1$, and the value of the
pullback of $f$ to $U\times_{\mathcal{M}}S''$ at $u''$ is again $\pm 1$. Finally, since
$u''$ is a point of total ramification for the morphism $S'' \to S'$, and the degree
of the morphism is $4$, the value of
$c(f)$ at $u'$ is $(\pm1)^4=1$. It follows that every function in $H^0(\mathcal{M};
\mathcal{I})$ takes the value $1$ at $u'$, and since these functions are constant, it
must be identically $1$. Thus $H^0(\mathcal{M}; \mathcal{I}) =\{1\}$ as claimed.
\end{proof}

\section{The general case}

Our goal in this section is to give a unified proof of Theorem~\ref{main}. Throughout
this section,  let 
$\rho\colon G
\to GL(V)$ denote the unique properly quasi-minuscule representation of
$G$. Hence
$\Ker
\rho=Z(G)$ and
$\overline{\mathcal{A}}_\rho = \mathcal{A}/Z(G)$. 

\subsection{Identification of the image torus}

Let $R_s$ be the set of short roots of $G$.  It is the set of  nonzero weights of
$\rho$. Let   $m_0 \geq 1$ be  the multiplicity of the trivial
weight in $\rho$. For  $\alpha\in R_s$, let $L_\alpha$ denote the
corresponding one-dimensional weight space of $V$. Let $L_0$ be the subspace of $V$
where $H$ acts trivially. 

\begin{lemma}\label{Hprimeforqm} $\rho(H)$ is  the subgroup $\widetilde{H}$ of 
$GL(V)$ consisting of the $A$ such that $A|L_0 =\Id$, 
$A|L_\alpha$ is multiplication by a nonzero scalar $a_\alpha$ and   the $a_\alpha$
satisfy: $a_{\alpha}a_{-\alpha} =1$, and, for every triple $(\alpha, \beta, \gamma)\in
R_s^3$ such that 
$\alpha+\beta+ \gamma =0$,
$a_{\alpha}a_{\beta} a_{\gamma}=1$.
\end{lemma}
\begin{proof} The element $\rho(h)$ acts on $L_0$ as the identity and on $L_\alpha$ via
$\alpha(h)$. Thus clearly $\rho(H)\subseteq \widetilde{H}$. To see the opposite
inclusion, we begin with the following general discussion.  Under the identification of
$W/W_0$ with
$R_s$, there is a homomorphism
$\varphi \colon \Zee[W/W_0] \to \Lambda^*$. Here, if $e_\alpha$ is the basis vector
corresponding to $\alpha$, $\varphi(e_\alpha) =\alpha$. The image of $\varphi$  is the
span of the short roots, which is just the root lattice $Q(R)$. The kernel
$K$ of
$\varphi$ is the set of all relations among the elements of $R_s$.

\begin{lemma} The kernel $K$ is generated by the following two types of relations:
\begin{enumerate}
\item[\rm  (i)] For $\alpha \in R_s$, $e_{\alpha} + e_{-\alpha}$;
\item[\rm  (ii)] For $\alpha, \beta, \gamma \in R_s$ such that $\alpha+\beta
=\gamma$,
$e_{\alpha} + e_{\beta}-e_{\gamma}$.
\end{enumerate}
\end{lemma}
\begin{proof}If $G$ is not simply laced, the set $R_s$ forms a root system whose span
is that of $R$. Thus it suffices to consider the case $R=R_s$, i.e.\  $G$ is simply
laced (but not necessarily irreducible).   Let
$v= \sum _{\alpha
\in R_s}n_\alpha e_\alpha
\in K$, so that $\sum _{\alpha \in R_s}n_\alpha \alpha =0$.  Using relations of type
(i), we may assume that all of the $\alpha$ are positive roots. If $\alpha_1, \dots,
\alpha_r$ are the simple roots, and $\alpha = \sum _in_i\alpha_i$ is a positive root,
let $\ell(\alpha) = \sum_in_i$. The proof is by induction on the largest value of
$\ell(\alpha)$ where $n_\alpha \neq 0$. If this number is greater than one, then
$\alpha = \alpha_i+\beta$ for some $i$, where $\ell(\beta) =\ell(\alpha)-1$. Applying
the relation $n_\alpha(e_\alpha - e_{\alpha_i} - e_\beta)$ then eliminates $e_\alpha$.
Eventually we reach the stage where all $\alpha$ are simple. Since the simple roots are
linearly independent, it follows that  at this stage $v=0$.
\end{proof}

Let $F$ be the free $\Zee$-module with basis $f_\alpha, \alpha\in R_s^+$ and
$f_{\alpha,\beta}$, where $\alpha, \beta$, and $\alpha + \beta \in R_s$. Let
$\psi\colon F \to \Zee[W/W_0]$ be defined by $\psi (f_\alpha) = e_{\alpha} +
e_{-\alpha}$ and $\psi(f_{\alpha,\beta}) = e_{\alpha} + e_{\beta}-e_{\alpha + \beta}$.
The preceding lemma implies that
$$F \xrightarrow{\psi} \Zee[W/W_0] \xrightarrow{\varphi} Q(R) \to 0$$
is exact, where $Q(R)$ is the root lattice.

Applying $\Hom(\cdot, \Cee^*)$ to this sequence gives an exact sequence of tori:
$$\{1\} \to \Hom(Q(R), \Cee^*) \to \Hom( \Zee[W/W_0], \Cee^*) \to \Hom(F, \Cee^*).$$
Moreover, $\Hom(Q(R), \Cee^*) =\overline{H}$ is the maximal torus in the adjoint form
$\rho(G)$ of $G$, and the induced homomorphism $H \to \overline{H} \to  \Hom(
\Zee[W/W_0], \Cee^*)$ is just the homomorphism which sends $h\in H$ to the homomorphism
$\Zee[W/W_0]\to \Cee^*$ defined by $\alpha \mapsto \alpha(h)$. It follows that, given
$a_\alpha\in \Cee^*$ for $\alpha \in R_s$ satisfying the conditions of the lemma, there
is an
$h\in H$ such that $\alpha(h) =a_\alpha$. Thus $\rho(h) = A$ as claimed.
\end{proof}

\subsection{The spectral cover and the basic correspondence}

Since $\rho$ is not minuscule, the spectral cover $T$ is of the form $m_0\cdot
\overline T_1 + T_0$ for a positive integer $m_0$. The divisor
$D=\overline T_1\cap T_0$ is nonempty, as is its preimage $\hat D$ in the normalization
$\hat T_0$ of $T_0$. As the weights of $\rho$ are invariant under $-1$, there is an
involution $\tau$ of $\hat T_0$.

We shall need the following lemma on roots:

\begin{lemma}\label{transitive} The Weyl group
$W$ acts transitively on the set of unordered triples $\{\alpha, \beta, \gamma\}$ with
$\alpha, \beta, \gamma\in R_s$ and $\alpha+ \beta+ \gamma=0$. The stabilizer in $W$ of
such a triple contains a cyclic permutation of the set $\{\alpha, \beta, \gamma\}$.
\end{lemma}

\begin{proof} To prove the first statement, it suffices to prove that $W$ acts
transitively on the set of unordered pairs
$\{\alpha, \beta\}$ such that  $\alpha$, $\beta$, and $\alpha+ \beta$ are short roots,
i.e.\ with  the set of embeddings of a root system of type $A_2$ into the subroot
system of short roots $G$ up to the action of the outer automorphism group of $A_2$.
Given two such pairs $\{\alpha_1, \beta_1\}$ and $\{\alpha_2, \beta_2\}$, we may assume
that $\alpha_1=\alpha_2$ since $W$ acts transitively on the set of short roots. Then
$\{\alpha_1, \beta_1, \beta_2\}$ spans a root system of type $A_2$ (if
$\beta_1=\beta_2$ or if the Cartan integer $n(\beta_1, \beta_2)=-1$) or of type $A_3$
(if the Cartan integer $n(\beta_1, \beta_2)$ is $0$ or $1$), and the result can be
checked directly in either case.

The second statement is an immediate consequence of the fact that the Weyl group of a
root system of type $A_2$ is the symmetric group $S_3$ and contains an index two
subgroup as claimed.
\end{proof}

We now define a correspondence on $\hat T_0$ which commutes with $\tau$ in an
appropriate sense. As we have seen in Lemma~\ref{transitive}, $W$ acts transitively on
the set of unordered triples $\{\alpha, \beta, \gamma\}$ with
$\alpha, \beta, \gamma\in R_s$ and $\alpha+ \beta+ \gamma=0$.  Let $\beta\in R_s$ be
such that $\varpi+\beta =-\gamma\in R_s$. (Such a $\beta$ exists except in case $G$ is
of type
$B_n$ or $A_1$.)  Let
$W'\subseteq W$ be the stabilizer of the set $\{\varpi, \beta, \gamma\}$, and set $S'
=\widetilde{\mathcal{M}}/W'$. Let  $\mathcal{R}_1$ be the subset of ordered triples
$(\alpha, \beta,
\gamma)\in R_s^3$ such that $\alpha+\beta+\gamma =0$.  Then $W$ acts diagonally on 
$\widetilde{M}\times \mathcal{R}_1$, and we let
$S'' =  (\widetilde{M}\times \mathcal{R}_1)/W$. One can check that $S''$ is irreducible,
i.e.\ that $W$ acts transitively on the set of ordered pairs $(\alpha,\beta)$ of short
roots whose sum is a short root, as long as $G$ is not of type $A_n$. However, we shall
not need this fact. Viewing $\hat T_0 =\widetilde{\mathcal{M}}/W_0$ as the quotient
$(\widetilde{\mathcal{M}}\times R_s)/W$ and   similarly for $S'$, we see that there are
morphisms $S''\to \hat T_0$ and $S''\to S'$. Thus as in the case of $E_6$ and $E_7$
there is a diagram
$$\begin{matrix}
& & S'' && &\\
& \swarrow && \searrow &\\
\hat T_0 & && & S'\\
& \searrow & &\swarrow &\\
& & \mathcal{M} && &
\end{matrix}.$$
There are obvious analogues of the involution $\tau$ on $S'$ and $S''$, for which the
above morphisms are equivariant. For example, the involution on $S''$ is induced by the
function
$(x, (\alpha, \beta, \gamma)) \mapsto (x,(-\alpha,
-\beta, -\gamma))$. We denote all of these involutions by $\tau$. Let $\nu'\colon S' \to
\mathcal{M}$ and $\nu''\colon S'' \to
\mathcal{M}$ be the induced morphisms, so that $\nu'\circ \tau = \nu'$, and similarly
for $\nu''$. We define the correspondence homomorphism $c\colon \hat\nu_*\mathbb{G}_{m,
\hat T_0} \to \nu_*'\mathbb{G}_{m,S'}$ in the usual way, by pulling a function up to
$S''$ and taking the image of the norm homomorphism from $\nu''_*\mathbb{G}_{m,S''}$ to
$\nu_*'\mathbb{G}_{m,S'}$. Clearly, $c$ commutes with $\tau$, and so there is an
induced homomorphism, which we also denote by $c$, from $(\hat\nu_*\mathbb{G}_{m,
\hat T_0})^- $ to $(\nu_*'\mathbb{G}_{m,S'})^-$.

For $G$ of type $A_1$ or $B_n$, there are no triples $\{\alpha, \beta, \gamma\}$ with
$\alpha, \beta, \gamma\in R_s$ such that $\alpha+ \beta+ \gamma=0$. In this case, the
constructions above are vacuous and there is no correspondence to consider. In this
case, we define $c$ to be the trivial homomorphism.

For future reference, we shall need the following lemma:

\begin{lemma}\label{subschemeX} Assume that $G$ is not of type $A_1$ or $B_n$. There is
a nonempty codimension two subscheme
$X$ of
$S'$ such that, for every \'etale open set $U\to \mathcal{M}$, if
$f\in \hat \nu_*\mathcal{B}_{\hat D}(U)$, then $c(f)|U\times _{\mathcal{M}}X=1$. If
$\beta \in R_s$ is   such that $\varpi+\beta =-\gamma\in R_s$   and
$W'$ is the stabilizer of the set $\{\varpi, \beta,\gamma\}$, then $X$ is
the subset of
$S'=\widetilde{\mathcal{M}}/W'$ where the functions from $S'$ to $E_{\rm reg}$
defined by $\varpi, \beta, \gamma$ are all equal to the identity.
\end{lemma}
\begin{proof} Let $\mathcal{R}_1$ be the subset of ordered triples $(\alpha, \beta,
\gamma)\in R_s^3$ such that $\alpha+\beta+\gamma =0$, and let $\mathcal{R}_2$ be the set
of unordered triples $\{\alpha, \beta, \gamma\}$ such that $\alpha,\beta, \gamma\in
R_s$ and $\alpha+\beta+\gamma =0$. There is the natural correspondence $\tilde c$ from
$\widetilde{\mathcal{M}}\times R_s$ to
$\widetilde{\mathcal{M}}\times \mathcal{R}_2$ defined as follows. Let $\varphi_1\colon
\widetilde{\mathcal{M}}\times \mathcal{R}_1\to \widetilde{\mathcal{M}}\times
R_s$ and $\varphi_2\colon \widetilde{\mathcal{M}}\times \mathcal{R}_1\to
\widetilde{\mathcal{M}}\times \mathcal{R}_2$ be the obvious morphisms: $\varphi_1(x,
(\alpha, \beta, \gamma)) = (x, \alpha)$ and $\varphi_2(x, (\alpha, \beta, \gamma)) = (x,
\{\alpha,
\beta, \gamma\})$. Then $\tilde c(g)$ is the pullback of $g$ under $\varphi_1$,
followed by the norm map induced by $\varphi_2$. Explicitly,
$$\tilde c(g)(x, \{\alpha, \beta, \gamma\}) = g(x, \alpha)g(x, \beta)g(x, \gamma).$$The
function
$g$ is the pullback of a function
$f$ on
$\hat T_0 = (\widetilde{\mathcal{M}}\times R_s)/W$, if and only if it is
$W$-invariant. In this case, 
$\tilde c(g)$ is also $W$-invariant, and hence is the pullback of a function on $S'$,
which is clearly $c(f)$.

Let $f$ be a function on $\hat T_0$ which is $1$ on $\hat D$. Then its pullback $g$ to
$\widetilde{\mathcal{M}}\times R_s$ is identically $1$ on $\coprod _{\alpha \in
R_s}D_\alpha$, where 
$$D_\alpha =\{ (x, \alpha) \in \widetilde{\mathcal{M}}\times R_s :
\alpha (x) \mbox{ is the identity $p_0$}\}.$$
Thus $\tilde c(g)$ is identically $1$ on the subscheme
$\coprod_{\{\alpha, \beta, \gamma\}\in \mathcal{R}_2}X_{\alpha, \beta, \gamma}$, 
where 
$$X_{\alpha, \beta, \gamma} = \{ (x, \{\alpha, \beta, \gamma\}) \in
\widetilde{\mathcal{M}}\times
\mathcal{R}_2 :
\alpha (x) = \beta(x) =\gamma(x)\mbox{ is the identity $p_0$}\}.$$
Clearly $X_{\alpha, \beta, \gamma}$ is nonempty and $W$-invariant, and it is of
codimension two since $\alpha+\beta+\gamma =0$. Hence $\coprod_{\{\alpha, \beta,
\gamma\}\in \mathcal{R}_2}X_{\alpha, \beta, \gamma}$ is the pullback of a nonempty
codimension two subset
$X$ of
$S'$ which has the required properties.
\end{proof}

We now claim:

\begin{proposition}\label{prop1} The sheaf $\overline{\mathcal{A}}_\rho$ is isomorphic to the
closed subsheaf $\Ker c \cap (\hat\nu_*{\mathcal B}_{\hat D})^-$ of $(\hat\nu_*{\mathcal
B}_{\hat D})^-$.
\end{proposition}
\begin{proof}
By Corollary~\ref{AbarinB},
$\overline{\mathcal{A}}_\rho$ is isomorphic to a closed subsheaf  of
$\hat\nu_*{\mathcal B}_{\hat D}$. Thus it suffices to prove that, via this
isomorphism,  
$\overline{\mathcal{A}}_\rho\cong \hat\nu_*{\mathcal B}_{\hat D}\cap (\hat\nu_*\mathbb{G}_{m,
\hat T_0})^-  \cap \Ker c$. Since both the image of $\overline{\mathcal{A}}_\rho$ and 
$\hat\nu_*{\mathcal B}_{\hat D}\cap (\hat\nu_*\mathbb{G}_{m,
\hat T_0})^-  \cap \Ker c$ are closed in $\hat\nu_*{\mathcal B}_{\hat D}$, it suffices
to prove that their restrictions to an open dense set agree. As usual, we use the open
set $U\to \mathcal{M}$ of Lemma~\ref{Gsplit} and Corollary~\ref{split}, and we shall
prove that their restrictions to $U$ agree. Note that $\hat \nu(\hat D)$ does not lie in
the image of $U$. Over $U$, $\mathcal{A}$ is the sheaf of sections of the
constant group scheme $U\times H$ and hence $\overline{\mathcal{A}}_\rho$ is the sheaf of
sections of the constant group scheme $U\times \overline{H}$. The
fiber of
$\hat T_0$ over a point
$x$ in the image of $U$ is indexed by the weights of $\rho$, i.e.\ by $R_s$. Thus if $f
\in
\hat\nu_*\mathbb{G}_{m, \hat T_0}(U)$, its value at $x$ is given by elements $a_\alpha
\in \Cee^*$. To say that $f\in  (\hat\nu_*\mathbb{G}_{m,
\hat T_0})^-  \cap \Ker c$ means exactly that $a_\alpha a_{-\alpha} =1$ and that, for 
every triple $(\alpha, \beta, \gamma)\in R_s^3$ such that 
$\alpha+\beta+ \gamma =0$,
$a_{\alpha}a_{\beta} a_{\gamma}=1$. By Lemma~\ref{Hprimeforqm}, this is exactly the
condition that $f$ lies in the image of $\overline{\mathcal{A}}_\rho$.
\end{proof}

\noindent \textbf{Proof of Theorem~\ref{main}.}  By Lemma~\ref{AvsAbar}, to prove
Theorem~\ref{main}, it will suffice to prove that
$H^1(\mathcal{M};\overline{\mathcal{A}}_\rho) =\{1\}$ and that
$H^0(\mathcal{M};\overline{\mathcal{A}}_\rho) =\{1\}$ or $\Zee$ depending on whether $E$ is
cuspidal or nodal. 
If $c$ is trivial, then $\overline{\mathcal{A}}_\rho =  (\hat\nu_*{\mathcal B}_{\hat
D})^-$ and the result is just Lemma~\ref{tau+1}. So we may assume that $c$ is not
trivial.  Restricting
$c$ to 
$(\hat\nu_*{\mathcal B}_{\hat D})^-$ and applying  Proposition~\ref{prop1} gives an
exact sequence
$$\{1\} \to \overline{\mathcal{A}}_\rho \to (\hat\nu_*{\mathcal B}_{\hat D})^- \to
\mathcal{I}
\to \{1\},$$
where $\mathcal{I}$ is the image of $(\hat\nu_*{\mathcal B}_{\hat D})^-$ in
$\nu'_*\mathbb{G}_{m, S'}$. 

As usual, assume first that $E$ is cuspidal. Since $S'$ is the quotient of the affine
space
$\frak h$ by a finite group, $H^0(\mathcal{M}; \nu'_*\mathbb{G}_{m, S'}) \cong \Cee^*$
is identified with the group of nonzero constant functions on $S'$. Since
$H^0(\mathcal{M};(\hat\nu_*{\mathcal B}_{\hat D})^- ) =
H^1(\mathcal{M};(\hat\nu_*{\mathcal B}_{\hat D})^- )=  \{1\}$ in this case, by
Lemma~\ref{tau+1}, it suffices to show that $H^0(\mathcal{M}; \mathcal{I}) = \{1\}$. In
any case, $H^0(\mathcal{M}; \mathcal{I})$ is a subsheaf of $H^0(\mathcal{M};
\nu'_*\mathbb{G}_{m, S'}) \cong \Cee^*$ and hence consists of constant functions  on
$S'$. To determine the value of these constant functions, it suffices to restrict to the
nonempty subset $X$ of Lemma~\ref{subschemeX}, where we have seen that the image of any
local section of
$\hat\nu_*{\mathcal B}_{\hat D}$ takes the value $1$. Hence $H^0(\mathcal{M};
\mathcal{I}) = \{1\}$, concluding the proof in the cuspidal case.

Finally  assume that  $E$ is nodal. In this case, by
Lemma~\ref{tau+1}, $H^0(\mathcal{M};(\hat\nu_*{\mathcal B}_{\hat D})^- ) \cong \Zee$
and 
$H^1(\mathcal{M};(\hat\nu_*{\mathcal B}_{\hat D})^- )=  \{1\}$.   We must again show
that
$H^0(\mathcal{M}; \mathcal{I})=\{1\}$. It suffices to show that
$$\{g\in H^0(S'; \mathbb{G}_{m, S'}): g|X =1\} =\{1\}.$$
We have seen that $S'=H/W'$, where $W'$ is the stabilizer of an unordered triple of
short roots
$\{\varpi, \beta,\gamma\}$ whose sum is zero. Thus 
$$H^0(S'; \mathbb{G}_{m, S'})
\subseteq (H^0(H; \mathbb{G}_{m, H}))^{W'} = (\Cee^*\times \mathsf{X}(H))^{W'} =
\Cee^*\times \mathsf{X}(H)^{W'}.$$
Suppose that $g\in H^0(S'; \mathbb{G}_{m, S'})$ and $g|X =1$. Then, since the origin
$1$ of $H$ is contained in the pullback of $X$, the pullback of
$g$ is the identity at $1\in H$. It follows that the pullback of $g$ lies in
$\mathsf{X}(H)^{W'}$. 

By Lemma~\ref{subschemeX}, the preimage of
$X$ in
$H$ contains the identity component of $\Ker \varpi \cap \Ker \beta \cap \Ker \gamma $,
which is a subtorus $H'$ of codimension two.  In particular, the pullback of $g$ is
equal to
$1$ on
$H'$. Let $H''$ be the subtorus of
$H$ whose Lie algebra is spanned by the coroots $\varpi\spcheck$, $\beta\spcheck$,  and
$\gamma\spcheck$. The induced homomorphism $H'\times H''\to H$ is a $W'$-equivariant
isogeny, so that 
$$\mathsf{X}(H)^{W'} \subseteq \mathsf{X}(H')^{W'}\oplus \mathsf{X}(H'')^{W'} .$$
Since $g$ has trivial restriction to $H'$, the image of $g$ lies in
$\mathsf{X}(H'')^{W'}$. Here $\mathsf{X}(H'') \cong \Zee^2$ is the weight lattice of a
root system of type $A_2$, and $W'$ contains a subgroup of order $3$ acting on $\Zee^2$
in the standard way. It follows that $\mathsf{X}(H'')^{W'} =0$, so that
$H^0(\mathcal{M}; \mathcal{I})=\{1\}$  as claimed.

This completes the proof of Theorem~\ref{main}.
\qed

\end{document}